\documentclass[12pt]{article}
\usepackage{amsmath,enumerate,amsfonts,color,amssymb, amsthm}

\usepackage{tikz}

\def\h{\mathcal{H}}
\def\NN{{\mathbb N}}
\def\ZZ{{\mathbb Z}}
\def\RR{{\mathbb R}}
\def\R{{\mathbb R}}

\def\CL{{\mycal{L}}}

\def\mcF{{\mycal F}}

\def\mcT{{\mycal T}}

\def\nd{\noindent}

\def\tiw{\tilde{w}}
\def\haw{\hat{w}}

\def\eps{{\varepsilon}}
\def\loc{{\rm loc}}

\newtheorem{theorem} {\sc  Theorem\rm} [section]

\newtheorem{corollary} [theorem] {\sc  Corollary\rm}
\newtheorem{lemma} [theorem] {\sc  Lemma\rm}
\newtheorem{proposition} [theorem] {\sc  Proposition\rm}

\newtheorem{definition}[theorem]{\sc  Definition\rm}
\newtheorem{open}[theorem]{\sc  Open problem\rm}
\newtheorem{remark}[theorem]{\sc  Remark\rm}

\def\nd{\noindent}

\newcounter{marnote}

\DeclareFontFamily{OT1}{rsfs}{}
\DeclareFontShape{OT1}{rsfs}{m}{n}{ <-7> rsfs5 <7-10> rsfs7 <10-> rsfs10}{}
\DeclareMathAlphabet{\mycal}{OT1}{rsfs}{m}{n}

\def\tr{{\rm tr}}

\def\mcS{{\mycal{S}}}

\def\Ss{{\mathbb S}}

\def\be{\begin{equation}}
\def\ee{\end{equation}}

\def\tr{{\rm tr}}

\def\mcS{{\mycal{S}}}

\def\mcE{{\mycal E}}

\def \f {\varphi}

\def\be{\begin{equation}}
\def\ee{\end{equation}}
\def\bea#1\eea{\begin{align}#1\end{align}}
\def\non{\nonumber}

\numberwithin{equation}{section}
\begin{document}

\title{Instability of point defects \\ in a two-dimensional nematic liquid crystal model\\
}
\author{Radu Ignat\thanks{Institut de Math\'ematiques de Toulouse, Universit\'e
Paul Sabatier, 31062
Toulouse, France. Email: Radu.Ignat@math.univ-toulouse.fr
}~, Luc Nguyen\thanks{Mathematical Insitute and St Edmund Hall, University of Oxford, Andrew Wiles Building, Radcliffe Observatory Quarter, Woodstock Road, Oxford OX2 6GG, United Kingdom. Email: luc.nguyen@maths.ox.ac.uk}~, Valeriy Slastikov\thanks{School of Mathematics, University of Bristol, University Walk, Bristol, BS8 1TW, United Kingdom. Email: Valeriy.Slastikov@bristol.ac.uk}~ and Arghir Zarnescu\thanks{University of Sussex, Department of Mathematics, Pevensey 2, Falmer, BN1 9QH, United Kingdom. Email: A.Zarnescu@sussex.ac.uk}\,\thanks{Institute of Mathematics ``Simion Stoilow" of the Romanian Academy, 21 Calea Grivitei Street, 01702 Bucharest, Romania}}

 \date{}

\maketitle
\begin{abstract}
We study a class of symmetric critical points in a variational $2D$ Landau - de Gennes model where the state of nematic liquid crystals is described by symmetric traceless $3\times 3$ matrices. These critical points play the role of topological point defects carrying a degree $\frac k 2$ for a nonzero integer $k$. We prove existence and study the qualitative behavior of these symmetric solutions. Our main result is the instability of critical points when $k\neq \pm 1, 0$.
\end{abstract}

%---------------------------------------------------%
%---------------------------------------------------%
\section{Introduction}
\label{sec:intro}

\subsection{Physical motivation}
\noindent The defining feature of  nematic liquid crystals is the local orientational ordering of the molecules. Its main macroscopic manifestation   is the emergence of certain patterns, called 
defects (points, lines or surfaces) where the local ordering, either  disappears or changes abruptly. Defects determine a number of the most important features of liquid crystals,  underlying spectacular phenomena and new prospective technologies, e.g.~knotted disinclination lines, bistable displays, control of nanoparticle suspensions (see \cite{PeterReview}).  These defects are often analysed in comparison with topological singular phenomena appearing in other fields of condensed and soft matter physics, such as superconductivity, materials science,  physics of polymers and even cosmology. 

There exist several competing  continuum liquid crystal theories describing the local orientational ordering by a specific {\it order parameter} (see \cite{dg,  ericksen, onsager-spatial_variation, of}). The most comprehensive and widely accepted continuum theory of nematic liquid crystals is the Landau-de Gennes theory \cite{dg}. It uses as an order parameter the so-called $Q$-tensor (a traceless, symmetric $3 \times 3$ matrix) so that the analysis is carried out in the five-dimensional space $\mcS_0$ of $Q$-tensors:
%and $L>0$ is an elastic constant (that from now on we will take for simplicity to be $1$). 
\begin{align*}
\mcS_0&=\bigg\{Q\in \R^{3\times3}\, :\,  Q=Q^t, \tr(Q)=0\bigg\}\\
&=\bigg\{s\left( n \otimes n -\frac{1}{3} I_3 \right)+r\left( m \otimes m -\frac{1}{3} I_3 \right)\, :\, s, r \in \R, n,m\in \Ss^2, n\cdot m=0\bigg\},
\end{align*}
where $\Ss^2$ is the unit sphere in $\R^3$, $I_3$ is the $3\times3$ identity matrix and $\big(n\otimes n \big)_{ij}=n_in_j$ for $1\leq i,j\leq 3$.

The central object in the Landau-de Gennes theory is the free energy functional $\mcF (Q)$; in fact, stable equilibrium configurations of the liquid crystalline system in $\Omega\subset\mathbb{R}^d$ ($d=2,3$) correspond to local minimisers of Landau-de Gennes energy.
In the simplest form, the free energy of a liquid crystal is given by
\begin{align}\label{LDG}
\mcF(Q)= \int_{\Omega} \Big[\frac{1}{2}|\nabla{Q}|^2 + f_{bulk}(Q)\Big]\,dx, \quad Q \in H^1_{loc}(\Omega, \mcS_0).
\end{align} 
The  simplest bulk potential $f_{bulk}(Q)$ that  captures the main physical characteristics is taken to be of the form 
$$
f_{bulk}(Q) = -\frac{a^2}{2}\tr(Q^2)-\frac{b^2}{3}\tr(Q^3)+\frac{c^2}{4}\left(\tr(Q^2) \right)^2,
$$
where $a^2, b^2, c^2>0$ are material constants. Note that the minimum set of the bulk potential $f_{bulk}(Q)$ is given by the set of uniaxial Q tensors (i.e., $Q$ has two equal eigenvalues):
\be \label{uniaxial}
\mcS_*=\left\{s_+\left(n\otimes n-\frac{1}{3}I_3\right)\, :\, n\in\mathbb{S}^2\right\}
\ee  with the constant order parameter $s_+$ given by
\be\label{def:s_+}
s_+=\frac{ b^2 + \sqrt{b^4+24 a^2 c^2}}{4 c^2}>0.
\ee

The critical points of the energy functional $ \mcF(Q)$ satisfy  the Euler-Lagrange equation:
\be\label{eq:EL}
 \Delta Q=-a^2 Q-b^2[Q^2-\frac 13|Q|^2I_3]+c^2Q|Q|^2 \quad \textrm{in }\, \Omega,
\ee where $\frac 13|Q|^2=\frac 1 3 \tr(Q^2)$ is the Lagrange multiplier associated to the traceless constraint. It is known that any $H^1_{loc}(\Omega, \mcS_0)$-solution of \eqref{eq:EL} is smooth in $\Omega$. The solutions of \eqref{eq:EL} describe the defects patterns, the simplest and most common being the point defects (see \cite{chandra, kleman, klelavr}). The analytical investigation of their structure and profile generates  very  challenging nonlinear analysis  problems.

The goal of this article is to investigate the profile and stability properties of point defects appearing for a certain type of symmetric solutions of \eqref{eq:EL} in the two dimensional case 
$$\Omega=\RR^2.$$ The boundary conditions imposed for these solutions are taken to be:
\be\label{BC1}
Q(x) \to  Q_k(x):=s_+ \left( n (x) \otimes n(x) -\frac{1}{3} I_3 \right) \quad \hbox{ as } |x| \to \infty,
\ee
where the map $n:\Omega\to \Ss^2$ is given in the polar coordinates by
\be\label{def:n}
n(x) =\left(\cos ({\textstyle\frac{k}{2}} \varphi) , \sin ({\textstyle\frac{k}{2}} \varphi) , 0\right), \quad r>0, \f\in [0, 2\pi),
\ee 
where $k \in \ZZ$ and $x=(r\cos \f, r\sin \f)$.
Note that if we consider $Q_k$ as an $\mathbb{R}P^1$-valued map on $\RR^2 \setminus \{0\}$, then it has degree $k/2$ about the origin. (For a definition of the degree for $\mathbb{R}P^1$-valued maps, see for instance \cite{BrezisCoronLieb}, p.$685-686$). This model can be seen as the $2D$ reduction of the physical situation of a $3D$ cylindrical boundary domain, with so-called ``homeotropic" boundary conditions
where the configurations are invariant in the vertical direction (see for instance \cite{bbch}).

\subsection{The $k$-radially symmetric solutions}

%{\red Radu: I don't like this phrase. Either you change it and make it clear, or otherwise please delete it : The most powerful feature of the $Q$-tensor description of liquid crystals is its ability of %describing much more complex patterns then possible in the other, simpler theories \cite{BallZar,baumanphilips, canevari, GolovatyMontero}.}
We will focus on the following type of symmetric solutions of \eqref{eq:EL} in the two-dimensional domain $\Omega=\RR^2$ that carry a topological information through the boundary condition \eqref{BC1}.  
%For all $k \in \ZZ \setminus \{ 0 \}$ and  $b^2 > 0$ we show existence of $k$-radially symmetric critical points of the Landau-de Gennes energy \eqref{LDG} in the whole space $\RR^2$ under boundary conditions \eqref{BC1}.  
\begin{definition}
\label{def:k_rad}
For $k \in \ZZ \setminus \{ 0 \}$, we say that a Lebesgue measurable map $Q:\Omega\to \mcS_0$ %(with $B_R(0)\subset\RR^2, R\in (0,\infty]$) 
is $k$-{\it radially symmetric} if the following conditions hold for almost every $x=(x_1, x_2)\in \Omega$:
\vskip 0.5cm

\noindent {\bf (H1)} The vector $e_3 =(0,0,1)$ is an eigenvector of $Q(x)$.

\noindent {\bf (H2)} The following identity holds 
$$
Q\bigg(P_2\big({\mathcal R}_k (\psi)  \tilde x\big)\bigg)= {\mathcal R}_k (\psi) Q(x) {\mathcal R}_k^t (\psi) , \ \textrm{for almost  every } \psi \in \RR, 
$$ where $\tilde x=(x_1,x_2,0)$, $P_2:\RR^3\to \RR^2$ is the projection given as $P_2(x_1,x_2,x_3)=(x_1,x_2)$ and
\begin{equation}
\label{eq: R_k }
{\mathcal R}_k (\psi) := \left(\begin{array}{ccc}\cos(\frac{k}{2}\psi) & -\sin(\frac{k}{2}\psi) & 0 \\\sin(\frac{k}{2}\psi) & \cos(\frac{k}{2}\psi) & 0 \\0 & 0 & 1\end{array}\right)
\end{equation} is the $\frac k 2$-winding rotation around the vertical axis $e_3$.
\end{definition}

\begin{remark}
If $k$ is an odd integer, then a map $Q \in H^1 (\Omega, \mcS_0)$ satisfying {\bf (H2)} automatically verifies {\bf (H1)} (see Proposition \ref{pro:H2}). 
%\begin{color}{red}luc: cut this out --Therefore,
%in this case, our results for $k$-{\it radially symmetric} maps $Q \in H^1 (B_R, \mcS_0)$ hold true by only assuming {\bf (H2)}. \end{color}
\end{remark}

\vskip 0.5cm We will show that the  $k$-radially symmetric solutions of \eqref{eq:EL} have a %deceptively simple-looking structure:
simple structure:

\begin{proposition} \label{2ode}
Let $k \in \ZZ \setminus \{ 0 \}$. If $Q\in H^1_{loc}(\RR^2, \mcS_0)$ is a $k$-radially symmetric solution of the Euler-Lagrange equations \eqref{eq:EL} on $\Omega=\R^2$ 
%(with $R\in (0,\infty]$ and  boundary conditions
%\be\label{BC}
%Q(x)=s_+ \left( n (\varphi) \otimes n(\varphi) -\frac{1}{3} I \right)\,\textrm{ for }|x|=R.
%\ee
satisfying the boundary conditions \eqref{BC1},
then $Q$ is smooth and  has the following form for every $x\in \RR^2$: 
\be\label{anY}
Q(x)= u(|x|) \sqrt{2}\left(n(x)\otimes n(x)-\frac{1}{2}I_2\right) + v(|x|) \sqrt{\frac{3}{2}}\left(e_3\otimes e_3-\frac{1}{3}I_3\right),
\ee
where $n$ is given in \eqref{def:n}, $I_2=I_3-e_3\otimes e_3$, $u \in C^2([0,\infty)) \cap C^\infty((0,\infty))$, $v \in C^\infty([0,\infty))$ and $u$ and $v$ satisfy the following system of ODEs in $(0, \infty)$: 
\be\label{ODEsystem}
\begin{cases}
u''+\frac{u'}{r}-\frac{k^2u}{r^2} &=u\left[-a^2+\sqrt{\frac{2}{3}} b^2 v+c^2\left( u^2+ v^2\right)\right]\\
v''+\frac{v'}{r}&=v\left[-a^2-\frac{1}{\sqrt{6}}b^2 v+c^2\left( u^2+ v^2\right) \right] + \frac{1}{\sqrt{6} } b^2 u^2,
\end{cases}
\ee
%Using boundary conditions at $R$, in order to obtain smooth solutions  and recalling that $Y$ has to solve \eqref{eq:EL}  we take the following boundary conditions on $(u,v)$
subject to boundary conditions:
\be\label{bdrycond}
u(0)=0,{\ v'(0)=0}, \ u(+\infty)=\frac{1}{\sqrt{2}} s_+,\,\,\,v(+\infty)=-\frac{1}{\sqrt{6}}s_+.
\ee
Conversely, if $u \in H_{loc}^1([0,\infty); rdr)  \cap L^2_{loc}([0,\infty);\frac{dr}{r})$ and $v \in H_{loc}^1([0,\infty);rdr)$ satisfy \eqref{ODEsystem} with the boundary condition $u(+\infty)=\frac{1}{\sqrt{2}} s_+$ and $v(+\infty)=-\frac{1}{\sqrt{6}}s_+$, then the tensor $Q$ defined by \eqref{anY} belongs to $H^1_{loc}(\RR^2,\mcS_0)$ and is a $k$-radially symmetric smooth solution of \eqref{eq:EL}-\eqref{BC1}.
\end{proposition}

Analysing the above ODE system, we construct solutions of \eqref{ODEsystem} - \eqref{bdrycond} using variational methods that lead to $k$-radially symmetric solutions of the 
Euler-Lagrange equations \eqref{eq:EL} with the boundary conditions \eqref{BC1}.

\begin{theorem}\label{thm:main0}
Let $a^2, b^2, c^2 >0$ be any fixed constants and $k \in \ZZ \setminus\{ 0 \}$. There exist $k$-radially symmetric solutions $Q \in H_{loc}^1(\RR^2, \mcS_0)$ of \eqref{eq:EL} - \eqref{BC1} having the form \eqref{anY}. Moreover, the corresponding profiles $(u,v)$ in \eqref{anY} satisfy the ODE system \eqref{ODEsystem} - \eqref{bdrycond} 
together with 
\[
u>0 \text{ and }v<0 \text{ in }(0,\infty).
\]
\end{theorem}

\begin{remark}
The case $b^2=0$ was studied in \cite{DRSZ} on bounded domains. They showed that on bounded domains, the ODE system has a unique solution under the assumption that $u > 0$ and $v < 0$. However, for infinite domain, the condition $b^2>0$ is essential in Theorem \ref{thm:main0}: there is \emph{no solution} to the ODE system on $(0,\infty)$ with $b^2 = 0$ which satisfies $u > 0$ and $v < 0$. See Appendix \ref{app:b0}.
\end{remark}

\begin{open}
For $b^2 > 0$, does the ODE system \eqref{ODEsystem}-\eqref{bdrycond} have a unique solution? (See Proposition \ref{pro:reg1} for a statement to this effect in a special case.)
\end{open}

\subsection{Instability of $k$-radially symmetric solutions for $k\in \ZZ \setminus \{ 0, \pm 1\}$}

Our main result concerns the local {instability} for {\bf all} $k$-radially symmetric critical points of $\mcF$ subject to \eqref{BC1} when $k \in \ZZ \setminus \{ 0, \pm 1\}$:
% and local stability of a special $k$-radially symmetric solution $Y$ defined in \eqref{anY} when $k = \pm 1$. 

\begin{theorem} \label{thm:main}
Let $a^2, b^2, c^2 >0$ be any fixed constants and $k \in \ZZ \setminus\{0,\pm 1 \}$. Any $k$-radially symmetric critical point $Q$ of \eqref{LDG} with $\Omega=\RR^2$ satisfying the boundary conditions \eqref{BC1} is locally unstable, i.e. there is a perturbation $P \in C_c^\infty(\RR^2, \mcS_0)$, supported in a bounded disk $B_R$, such that the second variation $\CL[Q] (P) <0$, where
\begin{align} \label{Eq:CLDef}
{\CL}[Q](P)&=\frac{1}{2}\frac{d^2}{dt^2}|_{t=0} \int_{\RR^2} \Big\{\frac{1}{2}|\nabla(Q+tP)|^2 + f_{bulk}(Q + tP) - \frac{1}{2}|\nabla Q|^2 - f_{bulk}(Q)\Big\}\,dx\non\\
&=\int_{\RR^2}\Big\{\frac{1}{2}|\nabla P|^2-\frac{a^2}{2}|P|^2-b^2\tr(P^2 Q)+\frac{c^2}{2}\left(|Q|^2|P|^2+2|\tr(QP)|^2\right)\Big\}\,dx.
\end{align}
%If $k=\pm 1$ then there exist a solution $(u,v)$ of \eqref{ODEsystem}, \eqref{bdrycond} such that $k$-radially symmetric $Q$-tensor $Y$ defined in \eqref{anY} is locally stable, i.e. the second variation $\CL[Y] (P) \geq 0$ for all $P \in H^1(\RR^2, \mcS_0)$. Moreover, $\CL[Y] (P) =0$ if and only if $P \in \{ \partial_{x_i} Y \}_{i=1}^2$, i.e. the kernel of the second variation coincides with translations of $Y$.
\end{theorem}

\begin{open}
Is it true that $k$-radially symmetric solutions of \eqref{eq:EL} in $\RR^2$ subject to \eqref{BC1} are stable for $k=\pm 1$?
\end{open}

\begin{remark}
%Due to the above theorem, for all $k \in \ZZ\setminus\{ 0, \pm 1 \}$, $b^2 > 0$ one expects a symmetry breaking for the ground state profile of the 2D point defect in $\RR^2$.  

This instability behaviour is drastically different from the case $b^2=0$ on a bounded disk $B_R$ centered at the origin. In \cite{DRSZ}, it was shown that the functional $\mcF$ with a boundary condition similar to \eqref{BC1} has a unique global minimizer in $H^1(B_R,\mcS_0)$, and furthermore that minimizer is $k$-radially symmetric. The deeper reason for this seems to be related to  the different structure of the minimum set of the potential $f_{bulk}$, which for $b^2=0$ is a $4D$ sphere while for $b^2\not=0$ is the $2D$ real projective plane.
\end{remark}

 There have been numerous numerical and analytical studies of two-dimensional point defects in the Landau - de Gennes framework \cite{baumanphilips, canevari, cladiskle, DRSZ, Fat-Slas, GolovatyMontero, Zhang, KraVirga, KraljVirgaZumer} {(also in micromagnetics see e.g. \cite{DIO, IgnOtt})}. Let us briefly mention a few papers that are directly related to this work. Our motivation came from the recent paper \cite{DRSZ}  where global minimisers of Landau-de Gennes energy are investigated on finite two-dimensional balls in the extreme low-temperature regime ($b^2=0$)  under $k$-radially symmetric homeothropic boundary conditions. The authors show that there exists a unique global minimizer of the Landau-de Gennes energy which is $k$-radially symmetric and provide the description of the ground state profile of a point defect of index $k/2$  in terms of the system of two ordinary differential equations (see \eqref{ODEsystem}). More general domains and boundary conditions were treated analytically (see \cite{baumanphilips, canevari, GolovatyMontero}) and numerically (see \cite{Zhang}). In \cite{baumanphilips} the Landau-de Gennes energy was investigated in a restricted three dimensional space of $Q$-tensors. The authors showed that in the case of small elastic constant the minimizers of Landau-de Gennes energy exhibit behavior similar to those of Ginzburg-Landau energy \cite{vortices}, namely for boundary conditions of degree $k/2$ there are exactly $k$ vortices of degree $\pm 1/2$. In \cite{canevari, GolovatyMontero} the minimizers of the full Landau-de Gennes energy were studied under non-orientable boundary conditions (which in our setting amounts to $k$ being odd). It was shown that in the low temperature regime and in the case of small elastic constant the minimizer has only one vortex.
 
The paper is organised as follows: in the next section we provide the basic properties of the $k$-radially symmetric maps that we study on balls $B_R$ of radius $R\in (0, \infty]$.  In Section $3$ we investigate the ODE system \eqref{ODEsystem} on bounded domains and prove certain fine qualitative properties of  solutions that will be used later. In Section $4$ we show the existence of a $k$-radially symmetric solution on the whole $\RR^2$ and investigate its behaviour at infinity. Finally, in Section 5 we investigate the stability of $k$-radially symmetric solutions and show Theorem~\ref{thm:main}. { Several open questions are also stated, some of them will be addressed in a forthcoming paper.}  

%In order to show the above results we 
%\begin{enumerate}
%\item need to prove existence of radial solution $(u,v)$ on the infinite domain
%\begin{itemize}
%\item Prove existence of a solution $(u,v)$ of E-L with $u>0$, $v<0$ for $R<\infty$ 
%\item Prove certain properties of $(u,v)$ for $R<\infty$. 
%\item Using upper and lower bounds on $(u,v)$ prove existence of solution for $R=\infty$.
%\item Show that for $R=\infty$ all solutions must have appropriate behaviour at $\infty$
%\end{itemize}
%\item using asymptotics at $\infty$ show that all radial solutions are unstable
%\end{enumerate}

\section{Basic aspects of $k$-radially symmetric maps, $k\neq 0$}

In order to classify $k$-radially symmetric maps on balls $B_R$ centered at the origin with $R\in (0, \infty]$ and $k\neq 0$ (see Definition~\ref{def:k_rad} for $\Omega=B_R$), we introduce some notation. We define $\{ e_i\}_{i=1}^3$ to be the standard basis in $\RR^3$ and denote, for $\varphi \in [0,2\pi)$,
\[
n=n(\f)=\left(\cos ({\textstyle\frac{k}{2}} \varphi) , \sin ({\textstyle\frac{k}{2}} \varphi) , 0\right), \,
m=m(\f)=\left(-\sin({\textstyle\frac{k}{2}} \varphi),\cos({\textstyle\frac{k}{2}} \varphi),0\right) .
\]
We endow the space $\mcS_0$ of $Q$-tensors with the scalar product $$Q\cdot \tilde Q=\tr(Q\tilde Q)$$ and 
for any $\f\in [0, 2\pi)$, we define the following orthonormal basis in $\mcS_0$:
\bea
E_0&=\sqrt{\frac{3}{2}}\left(e_3\otimes e_3-\frac{1}{3}I\right),\\
E_1&=E_1(\f)=\sqrt{2}\left(n\otimes n-\frac{1}{2}I_2\right),\,E_2=E_2(\f)=\frac{1}{\sqrt{2}}\left(n\otimes m+m\otimes n\right), \non \\
E_3&=\frac{1}{\sqrt{2}}(e_1 \otimes e_3+e_3\otimes e_1),\, E_4=\frac{1}{\sqrt{2}}\left(e_2\otimes e_3+e_3\otimes e_2\right).
\eea
Obviously, only $E_1$ and $E_2$ depend on $\f$ and we have 
\be
\label{deriv_E12}
\frac{\partial E_1}{\partial \f}=kE_2 \quad \textrm{and} \quad  \frac{\partial E_2}{\partial \f}=-kE_1.
\ee

We prove the following characterization of property {\bf (H2)} for $H^1(B_R, \mcS_0)$-maps.
\begin{proposition} \label{pro:H2}
Let $R\in (0, \infty)$, $k\neq 0$ and $Q \in  H^1(B_R, \mcS_0)$ be a map that satisfies {\bf (H2)} in $B_R$.Then $Q$ can be represented for a.e. $x=r(\cos \f, \sin \f)\in B_R$: 
\footnote{In these decompositions of $Q(x)$, the angle $\f$ defining $E_1$ and $E_2$ is given by the phase of $x\in B_R$.}
\begin{enumerate}
\item[1.] If $k$ is odd, then
$$
Q=\sum_{i=0}^2 w_i(r)E_i,
$$ 
where $w_0 \in H^1((0,R); r\,dr)$ and $w_1,w_2  \in H^1((0,R); r\,dr)\cap  L^2\left((0,R); \frac{1}{r}\,dr\right)$.

\item[2.] If $k$ is even, then 
$$
Q=\sum_{i=0}^2 w_i(r)E_i+(\tiw(r)\cos \frac k 2\f+\haw(r) \sin \frac k 2\f)E_3+(-\haw(r)\cos \frac k 2\f+\tiw(r) \sin \frac k 2\f)E_4,
$$ 
\end{enumerate}
where $w_0 \in H^1((0,R); r\,dr)$ and $\tiw, \haw, w_1,w_2  \in H^1((0,R); r\,dr)\cap  L^2\left((0,R); \frac{1}{r}\,dr\right)$. 
\end{proposition} 
\begin{proof}
Any $Q \in  H^1(B_R, \mcS_0)$ can be represented as
$$
Q(x)= \sum_{i=0}^4 w_i(x) E_i, \quad x\in B_R,
$$
with $w_i=Q\cdot E_i$ for $i=0,\dots, 4$ and
\begin{align}
|Q|^2&=\sum_{i=0}^4 w_i^2\non\\ \textrm{and}\quad |\nabla Q|^2&=\sum_{i=0}^4 |\nabla w_i|^2+\frac{k^2}{r^2}(w_1^2+w_2^2)+\frac{2k}{r^2}(
w_1 \frac{\partial w_2}{\partial \f}-w_2 \frac{\partial w_1}{\partial \f}),\label{Eq:Phantom}
\end{align}
where we used \eqref{deriv_E12}. 
Now we compute for $\psi\in \RR$:
\begin{align*}
&{\mathcal R}_k (\psi)e_1=n(\psi), {\mathcal R}_k (\psi)e_2=m(\psi), {\mathcal R}_k (\psi)e_3=e_3,  \\
&{\mathcal R}_k (\psi)n(\f)=n(\f+\psi), {\mathcal R}_k (\psi)m(\f)=m(\f+\psi), \f\in [0, 2\pi),
\end{align*}
so that we have for a.e. $x=r(\cos \f, \sin \f)\in B_R$:
\begin{align*}
{\mathcal R}_k(\psi) Q(x) {\mathcal R}_k^t(\psi) &=w_0(x)E_0+w_1(x)E_1(\f+\psi)+w_2(x)E_2(\f+\psi)\\
&+w_3(x)(\cos \frac k 2 \psi E_3+\sin \frac k 2 \psi E_4)+ w_4(x)(-\sin \frac k 2 \psi E_3+\cos \frac k 2 \psi E_4).
\end{align*}
Therefore, hypothesis {\bf (H2)} is equivalent (in polar coordinates) with:
\begin{align*}
& w_i(r, \f +\psi) = w_i (r, \f),  \  i =0, 1, 2\\
& w_3(r, \f +\psi) = w_3 (r, \f) \cos \frac k 2 \psi- w_4 (r, \f) \sin \frac k 2 \psi, \\
& w_4(r, \f +\psi) = w_3 (r, \f) \sin \frac k 2 \psi+ w_4 (r, \f) \cos \frac k 2 \psi, 
\end{align*}
for a.e. $r\in (0, R), \f \in (0, 2\pi), \psi\in \RR$. Therefore, we deduce that $w_i$ are independent of the angular variable $\f$ for $i=0, 1, 2$. Since $Q\in H^1(B_R, \mcS_0)$, we obtain that
$w_0 \in H^1((0,R); r\,dr)$, $w_i  \in H^1((0,R); r\,dr)\cap  L^2\left((0,R); \frac{1}{r}\,dr\right)$, $i=1,2$ and $w_3, w_4 \in H^1(B_R)$. It remains to characterize $w_3$ and $w_4$. Let $r\in (0,R)$ so that $w_3$ and $w_4$ are continuous on $\partial B_r$. (This is true because $w_3, w_4 \in H^1(\partial B_r)\subset C^{0, \frac 1 2}(\partial B_r)$ for a.e. $r\in (0,R)$.) Then the above equalities for $w_3$ and $w_4$ hold for every $\f \in [0, 2\pi)$ and $\psi \in \RR$. Setting $\tiw(r)=w_3(r, 0)$ and $\haw(r)=-w_4(r,0)$, we get that 
\begin{align*}
& w_3(r, \psi) = \tiw(r)  \cos \frac k 2 \psi+\haw(r) \sin \frac k 2 \psi, \\
& w_4(r, \psi) = \tiw(r) \sin \frac k 2 \psi-\haw(r) \cos \frac k 2 \psi, 
\end{align*}
for every $\psi \in \RR$. If $k$ is odd, the continuity of the $2\pi$-periodic functions $w_3(r, \cdot)$ and $w_4(r, \cdot)$ implies that $\tiw(r)=\haw(r)=0$ for a.e. $r\in (0,R)$. 
If $k$ is even, then 
\begin{align*}
&w_3^2(x)+w_4^2(x)=\tiw^2(r)+\haw^2(r)\\
& |\nabla w_3|^2(x)+|\nabla w_4|^2(x)=(\tiw')^2(r)+(\haw')^2(r)+\frac{k^2}{4r^2}(\tiw^2(r)+\haw^2(r))
\end{align*}
for a.e. $x=r(\cos \f, \sin \f)\in B_R$. The proof is now completed.
\end{proof}

As a consequence, we deduce the following characterization of $k$-radially symmetric maps defined on balls $B_R$:
\begin{corollary} \label{basis-id}
Let $R\in (0, \infty)$ and $k\neq 0$. A map $Q \in  H^1(B_R, \mcS_0)$ is $k$-radially symmetric if and only if $Q$ can be represented as 
\be \label{rel:basisdecomp}
Q=\sum_{i=0}^2 w_i(r)E_i, \quad x=r(\cos \f, \sin \f)\in B_R,
\ee 
where $w_0 \in H^1((0,R); r\,dr)$ and $w_i  \in H^1((0,R); r\,dr)\cap  L^2\left((0,R); \frac{1}{r}\,dr\right)$, $i=1,2$. Moreover, we have
\bea
\nonumber
\frac{\mcF(Q)}{2\pi}&=\int_0^R \left[ \frac{1}{2} \left( \sum_{i=0}^2 (w_i')^2+\frac{k^2}{r^2}(w_1^2+w_2^2) \right) -\frac{a^2}{2}\sum_{i=0}^2 w_i^2+\frac{c^2}{4}\left(\sum_{i=0}^2 w_i^2\right)^2\right]\,rdr\non\\
\nonumber
&-\frac{b^2}{3 \sqrt{6} }\int_0^R w_0\big[w_0^2-3(w_1^2+w_2^2)]\,rdr .
\eea
\end{corollary}

\begin{proof}
Assume that $Q$ is $k$-radially symmetric. If $k$ is odd, \eqref{rel:basisdecomp} follows directly from Proposition \ref{pro:H2}. If $k$ is even, by ${\bf (H1)}$, $e_3$ is an eigenvector of $Q$ and so the functions $\tilde w$ and $\hat w$ obtained in Proposition \ref{pro:H2} are zero almost everywhere in $(0, R)$. In either case, we have proved \eqref{rel:basisdecomp} . The converse implication is obvious. 
Now, we compute for $Q$ given by \eqref{rel:basisdecomp}:
$$\tr(Q^3)=\frac{1}{\sqrt{6} } w_0\big[w_0^2-3(w_1^2+w_2^2)].$$ The expression of $\mcF$ immediately follows.
\end{proof}

%Using Euler-Lagrange equations  \eqref{eq:EL}, \eqref{BC}  we obtain that $k$-radially symmetric critical points  correspond to solutions of the following system of ODEs (for odd $k$ we take $w_3=w_4=0$)
%\bea \label{5ode}
%& w_0 '' + \frac{w_0'}{r} = w_0 \left( -a^2+c^2\left(\sum_{i=0}^4 w_i^2\right) - \frac{b^2}{\sqrt{6}} w_0 \right) + \frac{b^2}{\sqrt{6}} (w_1^2 + w_2^2) - \frac{b^2}{2 \sqrt{6}} (w_3^2 + w_4^2) \\
%&w_1'' + \frac{w_1'}{r} -\frac{k^2}{r^2}w_1=w_1\left(-a^2+ c^2\left(\sum_{i=0}^4 w_i^2\right) + \frac{2 b^2}{\sqrt{6}} w_0 \right) - \frac{b^2}{2\sqrt{2}} (w_3^2 - w_4^2) \\
%&w_2'' + \frac{w_2'}{r}-\frac{k^2}{r^2}w_2= w_2\left( -a^2 + c^2\left(\sum_{i=0}^4 w_i^2\right) + \frac{2 b^2}{\sqrt{6}} w_0\right) - \frac{b^2}{\sqrt{2}} w_3 w_4\\
%&w_3'' + \frac{w_3'}{r}-\frac{k^2}{4 r^2}w_3=w_3 \left(-a^2+c^2\left(\sum_{i=0}^4 w_i^2\right)  - \frac{b^2}{\sqrt{6}} w_0   - \frac{b^2}{\sqrt{2}} w_1\right)   - \frac{b^2}{\sqrt{2}} w_2 w_4\\
%&w_4'' + \frac{w_4'}{r}-\frac{k^2}{4 r^2}w_4=w_4\left( -a^2+ c^2\left(\sum_{i=0}^4 w_i^2 \right) - \frac{b^2}{\sqrt{6}} w_0   + \frac{b^2}{\sqrt{2}} w_1\right)   - \frac{b^2}{\sqrt{2}} w_2 w_3
%\eea with boundary conditions  \eqref{BC}  becoming
%\be \label{5odeBC}
%w_0(R)=-\frac{s_+}{\sqrt{6}}, \ w_1(R)=\frac{s_+}{\sqrt{2}}, \ w_2(R)=w_3(R)=w_4(R)=0.
%\ee
We now provide the proof of  Proposition~\ref{2ode} with the characterization of $k$-radially symmetric solutions of \eqref{eq:EL}-\eqref{BC1}. In fact, 
we will prove the result for arbitrary balls $B_R$ with $R\in (0, \infty]$. The existence of such solutions is postponed to the next two sections.

\begin{proposition} \label{pro:sol_rad_sym} Let $k\neq 0$ and $R\in (0, \infty]$. If $Q\in H^1_\loc(B_R, \mcS_0)$ be a $k$-radially symmetric solution of the Euler-Lagrange equations \eqref{eq:EL} for $\Omega=B_R$
that satisfies the homeotropic boundary condition
\begin{equation}
Q(x)=s_+ \left( n (x) \otimes n(x) -\frac{1}{3} I_3 \right) \quad \textrm{on}\quad \partial B_R
	\label{BC1f}
\end{equation}
(with the convention \eqref{BC1} if $R=\infty$), then
$Q$ is smooth and
$$Q(x)=v(r) E_0+ u(r) E_1(\f) \quad \textrm{ for every } x=r(\cos \f, \sin \f)\in B_R,$$
$u \in C^2([0,R)) \cap C^\infty((0,R))$, $v \in C^\infty([0,R))$and the couple $(u,v)$ satisfies the ODE system \eqref{ODEsystem} and the boundary conditions
$u(0)=v'(0)=0$, $u(R)=\frac{1}{\sqrt 2}s_+$ and $v(R)=-\frac{1}{\sqrt 6}s_+$. 

Conversely, if $u \in H_{loc}^1([0,R); rdr)  \cap L^2_{loc}([0,R);\frac{dr}{r})$ and $v \in H_{loc}^1([0,R);rdr)$ satisfy \eqref{ODEsystem} with the boundary condition $u(R)=\frac{1}{\sqrt{2}} s_+$ and $v(R)=-\frac{1}{\sqrt{6}}s_+$, then the tensor $Q = v(r) E_0+ u(r) E_1(\f)$ belongs to $H^1_{loc}(B_R,\mcS_0)$ and is a $k$-radially symmetric solution of \eqref{eq:EL} and \eqref{BC1f}.
\end{proposition}

\begin{proof} 
Assume that $Q\in H^1_\loc(B_R, \mcS_0)$ be a $k$-radially symmetric solution of \eqref{eq:EL} and \eqref{BC1f}. Then $Q$ can be expressed in the form \eqref{rel:basisdecomp}. Standard elliptic regularity implies interior smoothness of any solution $Q\in H^1_\loc(B_R)$ of \eqref{eq:EL} (see for instance \cite{Ma-Za}). In particular, $w_i = Q \cdot E_i$ are smooth on $(0,R)$. We prove the remaining claim in several steps:

\medskip

\nd {\it Step 1: We prove that $w_1' w_2 -  w_2' w_1=0$ in $(0,R)$}. Fix $0<r_0<R$. Using \eqref{deriv_E12} and Corollary~\ref{basis-id} for $Q$ that is a (smooth) $k$-radially symmetric map in $B_{r_0}$, one computes that 
$$
\partial_\f Q = \frac{k}{2} \left( S_0 Q - Q S_0 \right) \quad \ \hbox{ in } B_{r_0},
$$
where
$
S_0= e_2 \otimes e_1 - e_1 \otimes e_2.
$
Considering now the scalar product of $S_0Q - QS_0$ with both parts of \eqref{eq:EL}, we obtain
$$
\Delta Q \cdot (S_0Q-QS_0) =0 \quad \textrm{in}\quad B_{r_0}.
$$
Integrating by parts over the ball $B_{r_0}$ leads to
$$
0=\int_{B_{r_0}} \Delta Q \cdot (S_0Q-QS_0) \, dx= \int_{\partial B_{r_0}} \partial_r Q \cdot (S_0Q -QS_0) \, d\h^1.
$$
Using the above expression of $\partial_\f Q$, we deduce
$$
\int_0^{2 \pi} \partial_r Q (r_0, \f) \cdot  \partial_\f Q (r_0, \f) \, d\f=0, \ \hbox{ for every } r_0\in(0,R).
$$
Combining with \eqref{deriv_E12} and \eqref{rel:basisdecomp}, we conclude with Step 1. 

\medskip

\nd {\it Step 2: We prove that $w_2=0$ in $(0, R)$.} First, note that the boundary conditions on $Q$ read as $w_0(R)=-\frac{s_+}{\sqrt{6}}$, $w_1(R) =\frac{s_+}{\sqrt{2}}$ and $w_2(R)=0$ (which are understood as limits if $R=\infty$). Next, we show that there exists $0<R_1<R$ such that $w_2(r)=0$ for all $r \in (R_1, R)$. Indeed, since $w_1$ is continuous and $w_1(R)>0$, there exists an interval $(R_1, R)$ such that $w_1>0$ on $(R_1, R)$. 
The equality in Step 1 implies that $\frac{w_2}{w_1}$ is constant on $(R_1, R)$ so that $w_2=0$ on $(R_1, R)$. 
In order to prove that $w_2$ vanishes in the whole interval $(0,R)$, we write the Euler-Lagrange equations \eqref{eq:EL} within the decomposition \eqref{rel:basisdecomp}:
\begin{align*}
& w_0'' + \frac{w_0'}{r} = w_0 ( -a^2+c^2 \sum_{i=0}^2 w_i^2 - \frac{b^2}{\sqrt{6}} w_0) + \frac{b^2}{\sqrt{6}} (w_1^2 + w_2^2), \\
&w_1'' + \frac{w_1'}{r} -\frac{k^2}{r^2}w_1=w_1\bigg(-a^2+ c^2\sum_{i=0}^2 w_i^2 + \frac{2 b^2}{\sqrt{6}} w_0\bigg), \\
&w_2'' + \frac{w_2'}{r}-\frac{k^2}{r^2}w_2= w_2\bigg( -a^2 + c^2\sum_{i=0}^2 w_i^2 + \frac{2 b^2}{\sqrt{6}} w_0\bigg)
\end{align*}
in $(0,R)$ where we used that $Q^2-\frac{|Q|^2}3 I_3=\frac{w_0^2-w_1^2-w_2^2}{\sqrt{6}}E_0-\sqrt{\frac 2 3}w_0 (w_1E_1+w_2E_2)$. Then
we apply the Cauchy-Lipschitz theory for the 2nd order ODE in $w_2$ (with smooth coefficients in $(0, R)$): since $w_2$ vanishes in $(R_1, R)$, we deduce that $w_2=0$ is the unique solution in $(0,R)$. Therefore, $Q=w_0(r)E_0+ w_1(r)E_1$ in $B_R$ and $(w_1, w_0)$ satisfies the system \eqref{ODEsystem}.

\medskip

\nd {\it Step 3: We prove $w_0'(0) = 0$ and the regularity of $w_0$.} Since $Q$ is smooth in $B_R$, we obtain that $w_0=Q\cdot E_0$ is smooth in $B_R$. In particular, $w_0$ extends to an even (smooth) function on $(-R,R)$. Therefore $w_0 \in C^\infty([0,R))$ and $w_0'(0) = 0$.

\medskip

\nd {\it Step 4: We prove that $w_1(0) = 0$ and the regularity of $w_1$.} By Corollary \ref{basis-id}, we know that $w_1\in H^1((0,R); r\,dr)\cap  L^2\left((0,R); \frac{1}{r}\,dr\right)$. Then $w_1$ is continuous on $(0,R)$ and we have for $r_1, r_2\in (0,R)$:
$$|w_1^2(r_2)-w_1^2(r_1)|=2\bigg|\int_{r_1}^{r_2}w_1w_1'\, dr\bigg|\leq 2\bigg(\int_{r_1}^{r_2}{w_1^2}\, \frac{dr}{r}\bigg)^{1/2} 
\bigg(\int_{r_1}^{r_2} (w'_1)^2\, rdr\bigg)^{1/2}.$$ Since the right hand side converges to zero as $|r_2-r_1|\to 0$, it follows that $w_1$ is continuous up to $r=0$. Combined again with $w_1\in L^2\left((0,R); \frac{1}{r}\,dr\right)$, we conclude that $w_1(0)=0$. 

For the regularity of $w_1$, note that $w_1$ satisfies 
\be \label{eq_interm}
w_1''+\frac{w_1'}{r}-\frac{k^2w_1}{r^2}={w_1(r)} g(r), \quad r\in (0,R)\ee  
where $g$ is a continuous function in $[0,R)$. Then we have (see \cite[Proposition 2.2]{ODE_INSZ}) that the function 
$$\zeta(r) = \frac{w_1(r)}{r^{|k|}}$$
 is continuously differentiable up to $r=0$ with vanishing derivative $\zeta'(0)=0$ at the origin. This implies that $\frac{w_1'}{r}-\frac{k^2w_1}{r^2}$ is continuous in $[0,R)$. Returning to equation \eqref{eq_interm}, we deduce that $w_1 \in C^2([0,R))$.
 
 Conversely, assume that $u \in H_{loc}^1([0,R); rdr)  \cap L^2_{loc}([0,R);\frac{dr}{r})$ and $v \in H_{loc}^1([0,R);rdr)$ satisfy \eqref{ODEsystem} with the boundary condition $u(R)=\frac{1}{\sqrt{2}} s_+$ and $v(R)=-\frac{1}{\sqrt{6}}s_+$. Then $Q$ belongs to $H^1_{loc}(B_R,\mcS_0)$ (thanks to \eqref{Eq:Phantom}) and satisfies \eqref{BC1f}. The system \eqref{ODEsystem} implies that $Q$ satisfies \eqref{eq:EL} in $B_R\setminus\{0\}$ (see Step 2 above). Since $Q \in H^1_{loc}(B_R,\mcS_0)$ and a point has zero Newtonian capacity in two dimensions, $Q$ satisfies \eqref{eq:EL} in $B_R$. This finishes the proof.
 \end{proof} 

\begin{proof}[Proof of Proposition~\ref{2ode}] It is a consequence of the above result.
\end{proof}

%---------------------------------------------------%
\section{Study of the ODE system on finite domains}
\label{sec:finiteODE}
In this section we first show the existence of a smooth solution $(u,v)$ of the system \eqref{ODEsystem} on a finite domain $(0,R)$ with
$k\in \ZZ\setminus\{0\}$ and $a^2, b^2, c^2>0$ with the boundary conditions
\be\label{BCodeFinite}
u(0)=0, v'(0)=0, u(R)=\frac{s_+}{\sqrt{2}}, v(R)= -\frac{s_+}{\sqrt{6}}. 
\ee
This solution $(u,v)$ has a sign invariance: $u>0$ and $v<0$ in $(0,R)$. Second, we study the qualitative properties and provide appropriate upper and lower bounds on the constructed solution $(u,v)$. These bounds will be extensively used in the next section when proving existence of the solution on the infinite domain. 

\subsection{Existence of solutions with $u>0$ and $v<0$}

Let $R\in (0, \infty)$, $k\neq 0$ and $a^2, b^2, c^2>0$. In order to prove existence of a solution $(u,v)$ of \eqref{ODEsystem} on $(0, R)$ satisfying \eqref{BCodeFinite} with the desired sign invariance, we will use a variational approach. First, note that a solution $(u, v)$ of the ODE system \eqref{ODEsystem} subject to $u(R)=\frac{s_+}{\sqrt{2}}$ and $v(R)= -\frac{s_+}{\sqrt{6}}$ is a critical point of the {\it reduced} energy functional: \footnote{If $k\neq 0$ and $R<\infty$, 
we have that $\mcE(u,v)<\infty$ if and only if $v\in H^1((0,R); r\,dr)$ and $u\in H^1((0,R); r\,dr)\cap  L^2\left((0,R); \frac{1}{r}\,dr\right)$. This is due to standard Sobolev embeddings and the fact that the bulk energy density is bounded from below (which can be seen from the inequality
$|v(v^2 -3u^2)|\leq (u^2+v^2)^{3/2}$ for any $u,v\in \RR$).\label{foot1}}
% and the existence of a constant $C_0>0$ such that $f(t)=-\frac{a^2}{2}t^2-\frac{b^2}{3 \sqrt{6} }t^3+\frac{c^2}{4}t^4\geq t^2-C_0$ for every $t=(u^2+v^2)^{1/2}\geq 0$. } 
\be \label{def:mcR}
\begin{split}
\mcE(u,v)&= \mcE_R(u,v) = \int_0^R \bigg[ \frac{1}{2} \left( (u')^2+(v')^2+\frac{k^2}{r^2}u^2 \right) -\frac{a^2}{2}(u^2 + v^2)+\frac{c^2}{4}\left(u^2+v^2\right)^2\\
& \quad -\frac{b^2}{3 \sqrt{6}}v(v^2 -3u^2)\bigg]\,rdr,
\end{split}
\ee
defined on the admissible set 
\be\label{SR}
\mcT= \left\{ (u,v) :  [0,R] \to \RR^2 \, \Big | \, \sqrt{r} u', \sqrt{r} v', \frac{u}{\sqrt{r}}, \sqrt{r} v \in L^2(0,R), \, u(R)=\frac{s_+}{\sqrt{2}}, v(R)= -\frac{s_+}{\sqrt{6}} \right\}.
\ee 
If $(u,v)\in \mcT$, then $u$ is continuous on $[0, R]$ with $u(0)=0$ (see Step 3 in the proof of Proposition \ref{pro:sol_rad_sym}).

%We now prove the existence of a couple $(u>0,v<0)$ that is solution of the ODE system \eqref{ODEsystem}:

\begin{proposition}\label{prop:positivlocalmin}
For any $0<R<\infty$ and $k\neq 0$, there exists a smooth local minimizer $(u,v)\in \mcT$ of $\mcE$ such that \eqref{BCodeFinite} holds
together with
 $u(r)>0$ on $(0,R]$ and $v(r)<0$ on $[0,R]$.  Moreover, $u\in C^2([0,R])$ with $\lim_{r\to 0} \frac{u}{r^{|k|}}>0$, $v\in C^2([0,R])$ and 
\be\label{rel:s3vr}
\sqrt{3}v(r)+u(r)<0, \forall r\in [0,R) .
\ee
\end{proposition}
\begin{proof} We divide the proof in several steps:

\medskip

\nd {\it Step 1: Existence of minimizers of $\mcE$ on 
\be
\label{T-}\mcT_-:=\{(u,v)\in \mcT \, :\, v \leq 0 \}.\ee} First, we know that $\mcE (u,v)$ is continuous and coercive in the convex closed set $\mcT_-$ endowed with the strong topology $\bigg(H^1((0,R); r\,dr)\cap  L^2\left((0,R); \frac{1}{r}\,dr\right)\bigg)\times H^1((0,R); r\,dr)$ (see Footnote \ref{foot1}). Then the direct method of calculus of variations implies the existence of a global minimizer $(u,v)$ of $\mcE $ on the subset $\mcT_-$. The couple $(u,v)$ satisfies 
\be\label{ODEineq}
\begin{split}
u''+\frac{u'}{r}-\frac{k^2u}{r^2} &={u}\left[-a^2+\sqrt{\frac{2}{3}} b^2 v+c^2\left( u^2+ v^2\right)\right],\\
v''+\frac{v'}{r}& {\geq} {v}\left[-a^2-\frac{1}{\sqrt{6}}b^2 v+c^2\left( u^2+ v^2\right) \right] + \frac{1}{\sqrt{6} } b^2 u^2
\end{split}
\, \textrm{ distributionally in } (0,R)
\ee
%Using boundary conditions at $R$, in order to obtain smooth solutions  and recalling that $Y$ has to solve \eqref{eq:EL}  we take the following boundary conditions on $(u,v)$
with boundary conditions 
\be\label{bdrycond-ineq}
u(0)=0, \ u(R)=\frac{1}{\sqrt{2}} s_+,\,\,\,v(R)=-\frac{1}{\sqrt{6}}s_+.
\ee
Since $u$ and $v$ are continuous in $(0,R]$, we have by \eqref{ODEineq} that $u\in C^2((0,R])$. 

Since the energy $\mcE$ is invariant with respect to a change of sign of $u$ and $u(R)>0$, we deduce that $(|u|,v)$ is also a global minimizer of $\mcE$ over $\mcT_-$. The strong maximum principle applied to the first equation in \eqref{ODEineq} (for $(|u|,v)$) implies 
$$|u| > 0 \quad \textrm{ in } \quad (0,R)$$
since $u(R)>0$. Hence
$$u> 0 \quad \textrm{ in } \quad (0,R).$$
Also note that on the open set $\{v<0\}$, the inequality in \eqref{ODEineq} becomes equality and therefore, $u, v\in C^\infty(\{v<0\}\cap (0,R))$.

\medskip

\nd {\it Step 2: We show that $\limsup_{r\to 0} v(r)<0$.} Assume by contradiction that $\limsup_{r\to 0} v(r)=0$. 
By construction, we know that $v$ is a global minimizer of $\mcE(u, \cdot)$ over the set $\mcT_-$. Note now that in $\mcE(u,v)$, the contribution of $v$ to the bulk potential is
\be \label{funct_f}
f^{(v)}(r, v)=r\bigg[\frac{b^2}{\sqrt{6}} u^2(r)v+\frac{-{a^2}+c^2u^2(r)}{2}v^2-\frac{b^2}{3 \sqrt{6} }v^3+\frac{c^2}4v^4\bigg], \quad r\in (0, R), v\leq 0.\ee Since $u(r)>0$ for all $r\in (0,R]$, we deduce the existence of $\delta\in(-\frac{1}{\sqrt{6}}s_+,0)$ such that $f(r,\cdot)$ is increasing in $[\delta, 0]$ for every $r\in (0,R)$.
(We highlight that $\delta$ depends only on $a^2, b^2, c^2>0$ and $ \| u \|_{L^\infty}$ and $\delta$ is independent in $r>0$ due to the form of the linear and quadratic terms in $f^{(v)}(r,v)$.) By the above assumption, there exists an interval $(R_1, R_2)\subset (0,R)$ such that $v(R_2)=\delta$, $v> \delta$ in $(R_1, R_2)$ and either $R_1=0$ or $v(R_1)=\delta$. Set $\tilde v=v$ in $(0,R_1)\cup (R_2, R)$ and $v=\delta$ in $(R_1, R_2)$. Then $\mcE(u, \tilde v)<\mcE(u, v)$ which contradicts the minimality of $v$.

\medskip

\nd {\it Step 3: We prove the following result: Let $(u,v)$ be a solution of \eqref{ODEineq} and \eqref{bdrycond-ineq} such that $u> 0$ and $v \leq 0$ in $(0,R)$. {Provided that $\limsup_{r\to 0} v(r)<0$}, then 
\eqref{rel:s3vr} holds true (which implies $v < 0$). Consequently, if $(u,v)$ is a minimizer of $\mcE$ in $\mcT_-$, then $(u,v)$ is a local minimizer of $\mcE$ in $\mcT$. } First, we define the function
\be
w=\frac{v}{u}+\frac{1}{\sqrt{3}} \quad \textrm{ in }\, (0,R).
\ee
Then one computes that
\be
w'=\frac{v'}{u}-\frac{u'v}{u^2}, \quad w''=\frac{v''}{u}-\frac{u''v}{u^2}-\frac{2u'v'}{u^2}+\frac{2(u')^2v}{u^3} \non
\ee 
that leads to
\be\label{eq:w0}
w''+\frac{w'}{r}+\frac{2u'}{u}w'=\frac{1}{u}\left(v''+\frac{v'}{r}\right)-\frac{v}{u^2}\left(u''+\frac{u'}{r}\right).
\ee
Using the ODE system \eqref{ODEineq} in  \eqref{eq:w0}, we obtain
\be\label{ineq:w}
w''+\left(\frac{1}{r}+\frac{2u'}{u}\right)w'-\underbrace{\frac{3b^2u}{\sqrt{6}}(-\frac{v}{u}+\frac{1}{\sqrt{3}})}_{\ge 0} w \geq\underbrace{-\frac{k^2v}{r^2u}}_{\ge 0} \quad \textrm{ in } (0,R).
\ee 
By definition of $w$, since $\limsup_{r\to 0} v(r)<0$, $u(0)=0$ and $u>0$ in $(0, R)$, we have $w<0$ in a neighborhood of $0$. By \eqref{bdrycond-ineq}, we also have $w(R)=0$. Applying the strong maximum principle for \eqref{eq:w0} on $(0,R)$, we deduce that $w < 0$ on $(0,R)$ and \eqref{rel:s3vr} is now proved.
 
%As a consequence, $v<0$ in $(0,R)$ since $u>0$ in $(0,R)$; thus, $(u,v)$ satisfies \eqref{ODEsystem}.

\medskip

\nd {\it Step 4: We prove the regularity of $u$, $v$ and $\lim_{r\to 0} \frac{u}{r^{|k|}}>0$.} By Proposition \ref{pro:sol_rad_sym}, the tensor $Q$ defined by \eqref{anY} is a smooth $k$-radially symmetric solution of \eqref{eq:EL}, and so by the same proposition, $u, v \in C^2([0,R])$ and $u(0) = v'(0) = 0$. Now, let $\zeta = \frac{u}{r^{|k|}}$, then $\zeta$ is continuous up to the origin (see the paragraph following \eqref{eq_interm}) and satisfies
\[
\zeta'' + (1 + 2|k|)\frac{\zeta'}{r} - g(r)\,\zeta(r) = 0 \text{ in } (0,R)
\]
for some function $g \in C([0,R))$. Applying \cite[Lemma B.2]{ODE_INSZ}, we see that $\zeta(0) > 0$, as desired.
\end{proof}

{
\begin{open}
By construction, the solution $(u,v)$ in Proposition \ref{prop:positivlocalmin} is a local minimizer of $\mcE$ over $\mcT$. Is $(u,v)$ a global minimizer?
\end{open}
}

\begin{remark} \label{pro:uvmax}
The following upper bound and uniqueness result are standard and hold for any $a^2, b^2, c^2>0$:
\begin{enumerate}
\item[1.] If $R\in (0, \infty)$ and $(u,v)$ is a solution of the ODE system \eqref{ODEsystem} subject to \eqref{BCodeFinite}, then the following upper bound holds: 
\be\label{bds:uvmax}
u^2+v^2\le\frac{2}{3}s_+^2 \quad \textrm{ in } (0,R);
\ee 
this remains true even if $b^2 = 0$ (see e.g. Proposition~3 in \cite{Ma-Za}).

\item[2.] There exists $R_0>0$ (depending on $a^2, b^2, c^2$) such that for any $R\in (0, R_0)$, there exists a unique solution $(u,v)$ of the ODE system \eqref{ODEsystem} with \eqref{BCodeFinite}. This is a consequence of the Poincar\'e inequality (see for instance in the related  Ginzburg-Landau framework Thm. $VIII.7$, p. $98$,\cite{vortices}).
\end{enumerate}
\end{remark}

\subsection{Upper and lower bounds for $(u,v)$}
Now we are ready to prove upper and lower bounds for \emph{any} solution $(u,v)$ of the ODE system \eqref{ODEsystem}-\eqref{BCodeFinite} with $u > 0$ and $v < 0$. These properties will be essential in proving the convergence of solutions on bounded domains to a solution on infinite domain. It turns out that these bounds strongly depend on the relation between material parameters $a^2$, $b^2$ and $c^2$.
In fact, we will distinguish regimes leading to different behavior of $v$ (see Figure \ref{figu}):

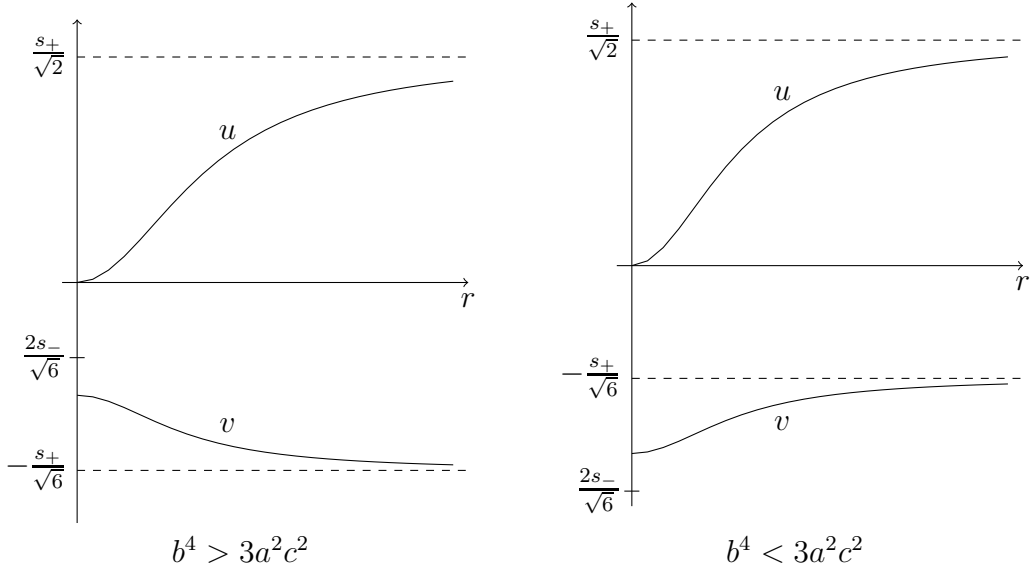
\begin{figure}
\begin{tabular}{ccc}

\begin{tikzpicture}
\draw[->] (-0.2,0)--(5.2,0) node[below] {$r$};
\draw[->] (0,-3.2)--(0,3.5);
\draw[dashed] (0,3)--(5.2,3);
\draw[dashed] (0,-2.5)--(5.2,-2.5);
\draw[domain=0:5,variable = \x] plot ({\x},{3*\x*\x/(\x*\x+3)});
\draw[domain=0:5,variable = \x] plot ({\x},{-1.5 - \x*\x/(\x*\x+2)});
\draw (0,-2.5) node[left] {$-\frac{s_+}{\sqrt{6}}$};
\draw (0,3) node[left] {$\frac{s_+}{\sqrt{2}}$};
\draw (-0.1,-1)--(0.1,-1);
\draw (0,-1) node[left] {$\frac{2s_-}{\sqrt{6}}$};
\draw (2,2) node {$u$};
\draw (2,-1.9) node {$v$};
\end{tikzpicture}

&\hspace{.5in}&
\begin{tikzpicture}
\draw[->] (-0.2,0)--(5.2,0) node[below] {$r$};
\draw[->] (0,-3.2)--(0,3.5);
\draw[dashed] (0,3)--(5.2,3);
\draw[dashed] (0,-1.5)--(5.2,-1.5);
\draw[domain=0:5,variable = \x] plot ({\x},{3*\x*\x/(\x*\x+2)});
\draw[domain=0:5,variable = \x] plot ({\x},{-2.5 + \x*\x/(\x*\x+2)});
\draw (0,-1.5) node[left] {$-\frac{s_+}{\sqrt{6}}$};
\draw (0,3) node[left] {$\frac{s_+}{\sqrt{2}}$};
\draw (-0.1,-3)--(0.1,-3);
\draw (0,-3) node[left] {$\frac{2s_-}{\sqrt{6}}$};
\draw (2,2.3) node {$u$};
\draw (2,-2.1) node {$v$};
\end{tikzpicture}\\
$b^4 > 3a^2c^2$ && $b^4 < 3a^2c^2$
\end{tabular}
\caption{Schematic graphs of $u$ and $v$ in different regimes of $a^2$, $b^2$ and $c^2$.}
\label{figu}
\end{figure}

\begin{itemize}
\item If $b^4  > 3a^2 c^2$ then $v \geq -\frac{s_+}{\sqrt{6}}$.
\item If $b^4 = 3 a^2 c^2$ then $v \equiv -\frac{s_+}{\sqrt{6}}$.
\item If $b^4 <  3a^2 c^2$ then $v \leq -\frac{s_+}{\sqrt{6}}$.
\end{itemize}

The regime $b^4=3a^2c^2$ can be considered as a special case the other regimes. However, it has a distinctive feature that $v= -\frac{s_+}{\sqrt{6}}$ is a local minimum of the $v$-relevant part of the bulk energy density (i.e. the function $f^{(v)}$ defined in \eqref{funct_f}). This allows us to establish stronger statement, for example the uniqueness result in Proposition \ref{pro:reg1} below.

 %{\red Radu: from here the physical interpretation should come...}

\subsubsection{The regime $b^4\geq3a^2c^2$}
Throughout this subsection we always assume 
\be\label{hyp:large}
b^4\geq3a^2c^2.
\ee
Under this assumption, the following inequalities hold (see \eqref{def:s_+})
\be\label{ineq:large}
\sqrt{\frac{2}{3}}s_-\geq -\frac{s_+}{\sqrt{6}} \geq -\frac{b^2}{\sqrt{6}c^2},
\ee where 
\be\label{def:s_-}
s_-=\frac{ b^2-\sqrt{b^4+24 a^2 c^2}}{4 c^2}<0.
\ee
When the inequality in \eqref{hyp:large} is strict, the inequalities in \eqref{ineq:large} are also strict.

We prove the following bounds on $u$ and $v$.
\begin{proposition}\label{lemma:propvforblarge}
Assume $b^4 \geq 3a^2c^2>0$, $0<R<\infty$, $k\neq 0$ and let $(u,v) \in \mcT$ be any solution of \eqref{ODEsystem}-\eqref{BCodeFinite} with $u > 0$ and $v < 0$ in $(0,R)$. Then 
\be\label{bds:vlarge}
-\frac{s_+}{\sqrt{6}}\le v\le \sqrt{\frac{2}{3}}s_-  \textrm{ in } (0,R)
\ee
and
\be\label{ineq-uu}
 u_{I}\leq u < \frac{s_+}{\sqrt{2}} \quad \hbox{ in } (0,R)
\ee 
where $u_{I}(r)$ is the unique solution of the following problem
\be \label{eq_u0}
u_{I}''+\frac{u_{I}'}{r}-\frac{k^2}{r^2}u_{I}=u_{I}\left[-a^2-\frac{b^2\sqrt{2}}{3}u_{I}+\frac{4c^2}{3}u_{I}^2\right], \quad u_{I}(0)=0,\ u_{I}(R) = \frac{s_+}{\sqrt{2}}.
\ee
\end{proposition}

\begin{proof} We divide the proof in several steps:

\medskip

\nd {\it Step 1: We prove the upper bound $v\le \sqrt{\frac{2}{3}}s_-$ in $(0,R)$.} Assume by contradiction that the maximum of $v$ is attained at some point $r_0\in [0,R)$ where  $0>v(r_0)> \sqrt{\frac{2}{3}}s_-$. By Proposition~\ref{prop:positivlocalmin}, we apply the maximal principle for the PDE satisfied by 
$v\in C^2(B_R)$:
\be \non
\underbrace{\Delta v(r_0)}_{\le 0}=\underbrace{v(r_0)\left[-a^2-\frac{b^2 v(r_0)}{\sqrt{6}}+c^2v^2(r_0)\right]}_{>0} +\underbrace{\left(v(r_0)c^2+\frac{b^2}{\sqrt{6}}\right)u^2(r_0)}_{>0}
\ee 
which leads to a contradiction.

\medskip

\nd {\it Step 2: We prove a weaker lower bound 
\be\label{est:vlowerbdsmall0}
v(r)\geq -\frac{b^2}{\sqrt{6}c^2} \textrm{ in } (0,R).
\ee
}
Assume by contradiction that the minimum of $v$ is achieved at $r_1\in [0,R)$ with $v(r_1)<-\frac{b^2}{\sqrt{6}c^2}$. Then, as at Step 1, we obtain 
\be \non
\underbrace{\Delta v(r_1)}_{\ge 0}=\underbrace{v(r_1)\left[-a^2-\frac{b^2 v(r_1)}{\sqrt{6}}+c^2v^2(r_1)\right]}_{<0} +\underbrace{\left(v(r_0)c^2+\frac{b^2}{\sqrt{6}}\right)u^2(r_0)}_{< 0}
\ee 
which leads to a contradiction. 

\medskip

\nd {\it Step 3: We prove the optimal lower bound $v(r)\geq -\frac{s_+}{\sqrt{6}}$ in $(0,R)$.}
Using \eqref{rel:s3vr} and  \eqref{est:vlowerbdsmall0} we obtain
\bea
\Delta v&=v\left[-a^2-\frac{b^2 v}{\sqrt{6}}+c^2v^2\right]+\left(vc^2+\frac{b^2}{\sqrt{6}}\right)u^2\non\\
&\le v\left[-a^2-\frac{b^2 v}{\sqrt{6}}+c^2v^2\right]+\left(vc^2+\frac{b^2}{\sqrt{6}}\right)3v^2\non\\
&\le v\left[-a^2+2\frac{b^2 v}{\sqrt{6}}+4c^2v^2\right] \quad \textrm{ in } B_R.\non
\eea 
Applying the maximum principle as at Step 2, we obtain the desired lower bound.

\medskip

\nd {\it Step 4: We prove $u(r) < \frac{s_+}{\sqrt{2}}$ in $(0,R)$.} Indeed, this upper bound follows directly from
inequalities \eqref{rel:s3vr} and \eqref{bds:vlarge}.

\medskip

\nd {\it Step 5: We prove the lower bound of $u$.} By \eqref{rel:s3vr} and \eqref{est:vlowerbdsmall0}, we have 
 $$
 v-\frac{1}{\sqrt{3}}u \geq 2v \geq -\frac{2 b^2}{\sqrt{6} c^2}  \quad \textrm{ in } (0,R).
 $$
 Multiplying with $v+\frac{u}{\sqrt{3}}<0$, we obtain: 
 $$
 \sqrt{\frac{2}{3}}b^2 v+c^2 v^2\le -\sqrt{\frac{2}{3}}b^2\frac{u}{\sqrt{3}}+\frac{c^2u^2}{3}  \quad \textrm{ in } (0,R).
 $$
By \eqref{ODEsystem}, the last inequality implies that $u$ is a super-solution for \eqref{eq_u0}, i.e.,
$$
u''+\frac{u'}{r}-\frac{k^2}{r^2}u\le u\left(-a^2-\frac{b^2\sqrt{2}}{3}u+\frac{4c^2}{3}u^2\right)  \quad \textrm{ in } (0,R).
$$ 
By \cite{Hervex2} (see also \cite{ODE_INSZ}), we know that there exists a unique solution $u_{I}$ of \eqref{eq_u0} that satisfies $0<u_{I}<\frac{s_+}{\sqrt{2}}$ in $(0,R)$; moreover, by \cite[Proposition 2.2]{ODE_INSZ}, we have that $u_{I}(r)=\alpha r^{|k|}+o(1)$ as $r\to 0$ for some $\alpha\geq 0$. By Step 4 in the proof of Proposition \ref{prop:positivlocalmin}, we deduce 
that $u(r)=\bar{\alpha} r^{|k|}+o(1)$ as $r\to 0$ with $\bar{\alpha}>0$, we can apply the comparison principle (see \cite[Proposition 3.5]{ODE_INSZ}) to obtain that
$u\geq u_{I}$ in $(0,R)$.
\end{proof}

When $b^4 = 3a^2c^2$, we have the following stronger result:
\begin{proposition} \label{pro:reg1}
Assume that  $b^4=3a^2c^2>0$, $k\neq 0$, $R\in (0, \infty)$ and let $(u,v) \in \mcT$ be any solution of \eqref{ODEsystem}-\eqref{BCodeFinite} with $u > 0$ and $v < 0$ in $(0,R)$. Then $v\equiv-\frac{s_+}{\sqrt{6}}$ and $u$ is 
the unique solution $u_{II}$ of the following problem:
\be \label{oldODE}
\begin{split}
&u_{II}''+\frac{u_{II}'}{r}-\frac{k^2u_{II}}{r^2}=c^2u_{II}(u_{II}^2-\frac{s_+^2}{2}), \quad \textrm{ in } (0,R),\\
& u_{II}(0)=0,\ u_{II}(R)=\frac{s_+}{\sqrt{2}}.
\end{split}
\ee
Moreover, $u_{II}$ is an increasing function with $0<u_{II} < \frac{s_+}{\sqrt{2}}$ on $(0,R)$.
\end{proposition}
\begin{proof} Note that $\sqrt{\frac{2}{3}}s_- = -\frac{s_+}{\sqrt{6}} = -\frac{b^2}{\sqrt{6}c^2}$. Thus, by Proposition \ref{lemma:propvforblarge} (namely \eqref{bds:vlarge}), $v \equiv -\frac{s_+}{\sqrt{6}}$. Substituting this in \eqref{ODEsystem}, we obtain that $u$ satisfies the problem \eqref{oldODE}.
By \cite{Hervex2} (see also \cite[Theorem~1.3]{ODE_INSZ}), we know that the problem \eqref{oldODE} admits a unique solution $u_{II}$. Moreover, $0<u_{II}<\frac{s_+}{\sqrt{2}}$ on $(0,R)$.
\end{proof}

\begin{remark}\label{rem:asym} In the case $R=\infty$, we recall that problem \eqref{oldODE} has a unique solution $u_{II}$ in $(0, \infty)$ 
(see \cite[Proposition~2.5]{ODE_INSZ}) and the behaviour of $u_{II}$ at infinity is given by:
$$
u_{II}(r) = \frac{s_+}{\sqrt{2}} - \frac{\beta}{r^2} + o\left(r^{-2}\right), \ \hbox{ as } r \to \infty,
$$
where
$\beta = \frac{k^2}{ \sqrt{2}b^2}$.
\end{remark}

%\begin{remark} \label{rem:bdsl}
%Let $(u,v)$ be the functions obtained in Proposition~\ref{prop:positivlocalmin}. Using Corollary~\ref{corollary:uuperbd},  Lemma~\ref{lemma:propvforblarge}, % Lemma~\ref{lemma:subsoluforblarge}  and \eqref{rel:s3vr} we obtain the following bounds on $(u,v)$
%\be
%u_0(r)\le u(r)\le\frac{s_+}{\sqrt{2}}, \quad  -\frac{s_+}{\sqrt{6}} \leq v\le -\frac{1}{\sqrt{3}}u_0.
%\ee
%\end{remark}

\subsubsection{The regime $b^4<3a^2c^2$ }
In this subsection we always assume 
\be\label{hyp:small}
b^4<3a^2c^2.
\ee
Under this assumption, the following inequalities hold (see \eqref{def:s_+}):
\be\label{ineq:small}
\sqrt{\frac{2}{3}}s_- < -\frac{s_+}{\sqrt{6}} < -\frac{b^2}{\sqrt{6}c^2}.
\ee
We prove the following bounds on $u$ and $v$.
\begin{proposition}\label{lemma:propvforbsmall}
Assume $0<b^4<3a^2c^2$, $0<R<\infty$, $k\neq 0$ and let $(u,v) \in \mcT$ be any solution of \eqref{ODEsystem}-\eqref{BCodeFinite} with $u > 0$ and $v < 0$ in $(0,R)$. Then  
\be\label{bounds:vsmall}
\sqrt{\frac{2}{3}}s_-\le v\le -\frac{s_+}{\sqrt{6}}
\ee
and
\be\label{subsol:usmall}
u_{III}\leq u<\frac{s_+}{\sqrt{2}} \quad \hbox{ in } (0,R).
\ee
where $u_{III}:(0,R)\to \RR$ is the unique solution of
\be
u_{III}''+\frac{u_{III}'}{r}-\frac{k^2u_{III}}{r^2}=\mu c^2 u_{III} (u_{III}^2 - \frac{s_+^2}{{2}}) , \quad u_{III}(0)=0,  \ u_{III}(R)=\frac{s_+}{\sqrt{2}}.
\ee with $\mu =  \frac{b^2}{\sqrt{b^4 + 24a^2c^2}}$.
\end{proposition}
\begin{proof}We follow the steps in the proof of Proposition \ref{lemma:propvforblarge}:

\medskip

\nd {\it Step 1: We prove the lower bound $v\geq \sqrt{\frac{2}{3}}s_-$ in $(0,R)$.} Assume by contradiction that the minimum of $v$ is achieved at some point $r_0\in [0,R)$ with $v(r_0)<\sqrt{\frac{2}{3}}s_-$. Then by \eqref{ineq:small}, the PDE satisfied by $v$ implies:
\be\non
\underbrace{\Delta v(r_0)}_{\ge 0}=\underbrace{v(r_0)[-a^2-\frac{b^2 v(r_0)}{\sqrt{6}}+c^2v^2(r_0)]}_{<0} +\underbrace{(v(r_0)c^2+\frac{b^2}{\sqrt{6}})u^2(r_0)}_{<0}
\ee 
which leads to a contradiction. 

\medskip

\nd {\it Step 2: We prove the weaker upper bound 
%\be\label{est:vupperblarge0}
$v(r)\le -\frac{b^2}{\sqrt{6}c^2}$ in $(0,R)$}. Assume by contradiction that the maximum of $v$ is achieved at $r_1\in [0,R)$ with $v(r_1)>-\frac{b^2}{\sqrt{6}c^2}$. Similarly, we obtain
\be\non
\underbrace{\Delta v(r_1)}_{\le 0}=\underbrace{v(r_1)[-a^2-\frac{b^2 v(r_1)}{\sqrt{6}}+c^2v^2(r_1)]}_{>0} +\underbrace{(v(r_0)c^2+\frac{b^2}{\sqrt{6}})u^2(r_0)}_{>0}
\ee 
which leads to a contradiction. 

\medskip

\nd {\it Step 3: We prove the optimal upper bound 
$v\le -\frac{s_+}{\sqrt{6}}$ in $(0,R)$.} By \eqref{rel:s3vr}, Step 2 leads to
\bea
\Delta v&=v[-a^2-\frac{b^2 v}{\sqrt{6}}+c^2v^2]+(vc^2+\frac{b^2}{\sqrt{6}})u^2\non\\
&\ge v[-a^2-\frac{b^2 v}{\sqrt{6}}+c^2v^2]+(vc^2+\frac{b^2}{\sqrt{6}})3v^2\non\\
&\ge v[-a^2+2\frac{b^2 v}{\sqrt{6}}+4c^2v^2].\non
\eea 
As above, the maximum principle yields the desired upper bound. 

%Using \eqref{bounds:vsmall} we obtain\marginnote{Arghir: this last inequality seems wrong}
%$$
%v''+\frac{v'}{r}\ge v[-a^2+\sqrt{\frac{2}{3}}b^2v+4c^2v^2]\ge 0
%$$ 
%and therefore a simple integration yields $v'\ge 0$ on $(0,R)$.

\medskip

\nd {\it Step 4: We prove the upper bound $u<\frac{s_+}{\sqrt{2}}$ in $(0,R)$}. The inequality $u\leq \frac{s_+}{\sqrt{2}}$ in $(0,R)$ follows directly from \eqref{bds:uvmax} and Step 3. Also, by \eqref{bounds:vsmall}, 
\[
-a^2 + \frac{2b^2}{\sqrt{6}}\,v + c^2\,v^2 \geq -\frac{1}{2}c^2\,s_+^2
\]
and so $u$ satisfies
\[
u'' + \frac{1}{r}u' - \frac{k^2}{r^2}\,u \geq c^2u(u^2 - \frac{1}{2}s_+^2).
\]
Thus, as $u \leq \frac{s_+}{\sqrt{2}}$ and $u(0) = 0$, the strong maximum principle implies that $u < \frac{s_+}{\sqrt{2}}$ in $(0,R)$.

\medskip

\nd {\it Step 5: We prove the lower bound of $u$.} First, note that \eqref{bds:uvmax} yields
\be \label{eq23}
-a^2+\frac{2b^2}{\sqrt{6}} v+c^2(u^2+v^2) \leq -a^2+\frac{2b^2}{\sqrt{6}} v+ \mu c^2(u^2+v^2) + (1-\mu)c^2 \frac{2}{3} s_+^2,
\ee
where $0<\mu<1$ will be chosen so that the function 
$$
\xi(v) = \frac{2b^2}{\sqrt{6}} v+ \mu c^2v^2
$$ 
is maximized on $[\sqrt{\frac{2}{3}}s_-, -\frac{s_+}{\sqrt{6}}]$ at the point $-\frac{s_+}{\sqrt{6}}$. For that, we need to insure that 
$\xi(-\frac{s_+}{\sqrt{6}}) > \xi( \sqrt{\frac{2}{3}} s_-)$ which is equivalent to
$$
\frac{2b^2}{\sqrt{6}} + \mu c^2 (  \sqrt{\frac{2}{3}} s_-  - \frac{s_+}{\sqrt{6}}) >0 \quad\textrm{i.e., } \quad (8+\mu)b^2 - 3\mu \sqrt{b^4+24 a^2c^2} >0
$$
(here we used \eqref{def:s_+} and \eqref{def:s_-}). Thus, the choice $\mu = \frac{b^2}{\sqrt{b^4 + 24a^2c^2}}\in (0,1)$ fulfills our objective and we conclude that
$$
\xi(v(r)) \leq \xi(-\frac{s_+}{\sqrt{6}}) = -\frac{b^2}{3} s_+ + \frac{\mu c^2}{6} s_+^2, \quad r\in (0,R).
$$
Combined with \eqref{eq23}, we obtain:
\be
-a^2+\frac{2b^2}{\sqrt{6}} v+c^2(u^2+v^2) \leq -a^2+ \mu c^2 u^2  -\frac{b^2}{3} s_+ + \frac{\mu c^2}{6} s_+^2+ (1-\mu)c^2 \frac{2}{3} s_+^2.
\ee
Together with \eqref{def:s_+}, it yields
%$$
%-a^2 -\frac{b^2}{3} s_+  + \frac{2c^2}{3} s_+^2 =0 \quad \hbox{and} \quad \frac{\mu c^2}{6} s_+^2 - \frac{2\mu c^2}{3} s_+^2 = -\frac{\mu c^2}{2} s_+^2.
%$$
%Therefore we have
$$
-a^2+\frac{2b^2}{\sqrt{6}} v+c^2(u^2+v^2) \leq \mu c^2 (u^2 - \frac{s_+^2}{2})
$$
and by \eqref{ODEineq},
\be
u''+\frac{u'}{r}-\frac{k^2}{r^2}u \leq \mu c^2 u ( u^2 - \frac{s_+^2}{{2}}), \quad r\in (0,R).
\ee 
As at Step 5 in Proposition \ref{lemma:propvforblarge}, we conclude to the desired lower bound using the comparison principle (see \cite[Proposition 3.5]{ODE_INSZ}). 
\end{proof}

%\begin{remark} \label{rem:bdss}
%Let $(u,v)$ be the ones obtain in  Proposition~\ref{prop:positivlocalmin}. Using   Lemma~\ref{lemma:propvforbsmall}, Lemma~\ref{lemma:subsoluforbsmall}  and 
%\eqref{bds:uvmax} we obtain the following bounds on $(u,v)$
%\be
%u_1(r)\le u(r)\le\frac{s_+}{\sqrt{2}}, \quad  -\sqrt{\frac{2}{3} s_+^2 -u_1^2(r)} \leq v (r) \le -\frac{s_+}{\sqrt{6}} .
%\ee
%\end{remark}

%---------------------------------------------------%
\section{Study of the ODE system on the infinite domain}
\label{sec:infinite}
In this section we study the ODE system \eqref{ODEsystem} on the infinite domain $(0,\infty)$ for $k\in \ZZ\setminus \{0\}$. Using results of the previous section, we first prove the existence of a solution of \eqref{ODEsystem} subject to \eqref{bdrycond}. As consequence, we prove existence of $k$-radially symmetric solutions of \eqref{eq:EL} on the whole $\RR^2$ stated in Theorem \ref{thm:main0}. Next, we prove finer asymptotic behavior at infinity of {\bf any} solution of \eqref{ODEsystem} subject to \eqref{bdrycond}.

\bigskip

We start by proving the following existence result on $(0,\infty)$.

\begin{proposition} \label{prop:inf}
Let $a^2, b^2, c^2 >0$ be fixed constants and $k \in \ZZ \setminus \{0\}$. Then there exists a smooth solution\footnote{Here, $u$ and $v$ are $C^2$ up to $r=0$.} $(u,v)$ of \eqref{ODEsystem} defined on $(0,\infty)$ with boundary conditions \eqref{bdrycond}. Moreover, we have $0<u< \frac{s_+}{\sqrt{2}}$, $v<0$ in $(0,\infty)$ and $(u,v)$ is locally minimizing in the following sense:
\[
\mcE_R(u,v) \leq \mcE_R(u + \xi, v + \eta) \text{ for any } (\xi, \eta) \in C_c^\infty(0,R) \text{ with } \sup_{(0,R)}|\eta| < \min\Big(\frac{s_+}{ \sqrt6}, \sqrt{\frac{2}{3}}|s_-|\Big),
\]
for any $R>0$, where $\mcE_R$ is given by \eqref{def:mcR}.

\end{proposition}
\begin{proof}
For every $n\in \NN^*$, let $(u_n, v_n)$ be the solution of \eqref{ODEsystem} on the interval $(0,n)$ subject to \eqref{BCodeFinite} constructed in Proposition~\ref{prop:positivlocalmin}. We extend $u_n$ and  $v_n$ to the infinite domain $(0, \infty)$  by setting the functions $(\bar u_n, \bar v_n):(0, \infty)\to \RR^2$ as follows:
$$
\bar u_n = 
\begin{cases}
u_n (r) & r \in (0,n) \\
\frac{s_+}{\sqrt{2}} & r  \geq n
\end{cases}, \quad 
\bar v_n = 
\begin{cases}
v_n (r) & r \in (0,n) \\
-\frac{s_+}{\sqrt{6}} & r  \geq n 
\end{cases}.
$$
Since $\{(\bar u_n, \bar v_n)\}_{n\geq 1}$ are uniformly bounded in $L^\infty (0, \infty)$, we have by standard regularity arguments that for any given compact interval $J \subset (0,\infty)$ and for large enough $n\geq n_J$, the couples $\{(\bar u_n,\bar v_n)\}_{n\geq n_J}$ are uniformly bounded in $C^3 (J)$. Using the Arzela-Ascoli theorem, we deduce that $\bar u_n \to u$ and $\bar v_n \to v$ in $C^2_{loc} (0,\infty)$ (up to a subsequence). Thus, $(u, v):(0, \infty)\to \RR^2$ satisfy \eqref{ODEsystem} on $(0, \infty)$, too. By Propositions \ref{pro:reg1}, \ref{lemma:propvforblarge}, \ref{lemma:propvforbsmall}, Step 3 in the proof of Proposition \ref{prop:positivlocalmin} and Remark \ref{pro:uvmax}, we have\footnote{For the strict upper bound of $u$, see Steps 4 in the proofs of Propositions \ref{lemma:propvforblarge} and \ref{lemma:propvforbsmall}.}
\begin{align*}
0 &\leq u < \frac{s_+}{\sqrt{6}},\\
\min\Big(-\frac{s_+}{ \sqrt6}, \sqrt{\frac{2}{3}}s_-\Big) &\leq v \leq \max\Big(-\frac{s_+}{ \sqrt6}, \sqrt{\frac{2}{3}}s_-\Big),\\
u + \sqrt{3}v & < 0,\\
u^2 + v^2 &\leq \frac{2}{3}s_+^2.
\end{align*}
Also, $u > 0$ in $(0,\infty)$.

%\marginnote{inserted}

We next show that $u \in H^1_{loc}([0,\infty); r\,dr) \cap L^2([0,\infty);\frac{dr}{r})$ and $v \in H^1_{loc}([0,\infty);r\,dr)$. Thanks to the (uniform) bound of $u_n$ and $v_n$, it suffices to show that $\mcE_{m}(u_n,v_n)$ is uniformly bounded for  $n > m \geq 0$. Indeed, if we $(\bar u_{m,n}(r),\bar v_{m,n}(r))$ be the extension of $(u_m,v_m)$ which equals to $(u_n,v_n)$ in the interval $(m+1,n)$ and is linear in $[m,m+1]$, then
\begin{align*}
0
	& \geq \mcE_n(u_n,v_n) - \mcE_n(\bar u_{m,n},\bar v_{m,n})\\
	&= \mcE_{m}(u_n,v_n) - \mcE_m(u_{m},v_{m})\\
		&\qquad + \int_{m}^{m+1} \bigg[{ \frac 1 2 (|u_n'|^2 + |v_n'|^2 + \frac{k^2}{r^2}|u_n|^2)} + f_{bulk}(u_n,v_n)\bigg]\,rdr \\
		&\qquad - \int_{m}^{m+1} \bigg[{ \frac1 2( |\bar u_{m,n}'|^2 + |\bar v_{m,n}'|^2 + \frac{k^2}{r^2}|\bar u_{m,n}|^2)} + f_{bulk}(\bar u_{m,n},\bar v_{m,n})\bigg]\,rdr,
\end{align*}
where (by a slight abuse of notation)
\[
f_{bulk}(x,y) = f_{bulk}(xE_1 + yE_0) = -\frac{a^2}{2}(x^2 + y^2)+\frac{c^2}{4}\left(x^2+y^2\right)^2 -\frac{b^2}{3 \sqrt{6}}y(y^2 -3x^2).
\]
As $(u_n,v_n)$ are uniformly bounded in $(m,m+1)$, $(\bar u_{m,n}(r),\bar v_{m,n}(r))$ and its derivative are also uniformly bounded in $(m,m+1)$. It thus follows that $0 \geq \mcE_{m}(u_n,v_n) - \mcE_m(u_{m},v_{m}) - C$ for some constant $C$ independent of $n$. This proves that $\mcE_{m}(u_n,v_n)$ is uniformly bounded.

The locally minimizing property of $(u,v)$ follows from the bounds for $v_n$ and the minimizing property of $(u_n,v_n)$. It remains to show that $(u,v)$ takes on the desired value at infinity (the boundary condition at the origin and the smoothness of $u$ and $v$ are a consequence of Proposition \ref{pro:sol_rad_sym}).

\medskip

\nd {\it Case 1: $b^4=3a^2c^2>0$.} By Proposition \ref{pro:reg1}, we have that $v=-\frac{s_+}{\sqrt{6}}$. Moreover, since $u_n$ is the unique solution of \eqref{oldODE} in $(0,n)$, we know by \cite[Proposition 2.4]{ODE_INSZ} that $u_n$ converges in $C_{loc}^2$ to the unique solution $u=u_{II}$ of \eqref{oldODE} in $(0, \infty)$, and so $u(\infty) = \frac{s_+}{\sqrt{2}}$.

\medskip

\nd {\it Case 2: $b^4>3a^2c^2>0$.} By Proposition \ref{lemma:propvforblarge}, the same argument as above implies that $-\frac{s_+}{\sqrt{6}}\le v\le \sqrt{\frac{2}{3}}s_-<0$ and $0<u_{I}\leq u \leq \frac{s_+}{\sqrt{2}}$ in $(0, \infty)$ where $u_{I}$ is the unique solution of \eqref{eq_u0} in $(0,\infty)$, in particular, $u(\infty)=\frac{s_+}{\sqrt{2}}$. For $v(\infty)$, note on one hand that $v\geq -\frac{s_+}{\sqrt{6}}$ in $(0, \infty)$ so that $\liminf_{r\to \infty} v\geq -\frac{s_+}{\sqrt{6}}$. 
On the other hand, by \eqref{rel:s3vr}, we have that $\limsup_{r\to \infty} v\leq-\frac1{\sqrt{3}}\lim_{r\to \infty} u=-\frac{s_+}{\sqrt{6}}$. We thus have $v(\infty) = -\frac{s_+}{\sqrt{6}}$.

%As at Step {\red 4} in Proposition \ref{prop:positivlocalmin}, we have that the function ${r}\mapsto \frac{u}{r^{|k|}}$ is continuously differentiable up to $r=0$, in particular $u$ is Lipschitz up to $r=0$ and $u(0)=0$. Moreover, 
%\eqref{rel:s3vr} and the bounds on $v$, imply that $u < \frac{s_+}{\sqrt{2}}$ in $(0, \infty)$. Furthermore, the argument in  Step~{\red 3} in the proof of Proposition \ref{pro:sol_rad_sym} implies that $v'(0)=0$. 

\medskip

\nd {\it Case 3: $0<b^4<3a^2c^2$.} Arguing as in the previous case (but using Proposition \ref{lemma:propvforbsmall} instead of Proposition \ref{lemma:propvforblarge}), we get $u(\infty) = \frac{s_+}{\sqrt{2}}$. Next, since $v\leq -\frac{s_+}{\sqrt{6}}$ in $(0, \infty)$, $\limsup_{r\to \infty} v \leq -\frac{s_+}{\sqrt{6}}$. On the other hand, by \eqref{bds:uvmax},
\[
\liminf_{r\to \infty} v \geq -\sqrt{\frac{2}{3}s_+^2 - u(\infty)^2} =  -\frac{s_+}{\sqrt{6}}.
\]
We again obtain $v(\infty)=-\frac{s_+}{\sqrt{6}}$ as desired.
\end{proof} 

We now prove the existence of $k$-radially symmetric solutions of \eqref{eq:EL} subject to \eqref{BC1}:

\begin{proof}[Proof of Theorem \ref{thm:main0}] 
The assertion is a consequence of Propositions \ref{2ode} and \ref{prop:inf}.
\end{proof}

% an
%Let $(u,v)$ be the solution of \eqref{ODEsystem} defined on $(0,\infty)$ subject to \eqref{bdrycond} that we have constructed at Proposition \ref{prop:inf}. Consider the map $Q:\RR^2\to \mcS_0$ defined by \eqref{anY}. Since $u$ and $v$ are Lipschitz continuous at $r=0$, we get by \eqref{bdrycond} that $Q\in H^1_{loc}(\RR^2, \mcS_0)$ and $Q$ satisfies \eqref{BC1}. By Corollary \ref{basis-id}, we deduce that $Q$ is a $k$-radially symmetric map. The same computation as Step 2 in Proposition \ref{pro:sol_rad_sym} shows that $Q$ is a solution of \eqref{eq:EL}. 

In the proof of the instability result, we need some detailed behavior at $\infty$ of {\bf any} solution $(u,v)$ of the system \eqref{ODEsystem} subject to \eqref{bdrycond}:

\begin{lemma} \label{lemma:asym}
Let $u$ and $v$ be {\bf any} solution
%\footnote{Here, we assume that $u$ is Lipschitz continuous and $v$ is $C^1$ up to $r=0$.} 
of \eqref{ODEsystem} defined on $(0,\infty)$ subject to \eqref{bdrycond}.
Then $(u,v)$ has the following behavior as $r \rightarrow \infty$:
\begin{align} \label{ab1}
u &= \frac{s_+}{\sqrt{2}} - \frac{\sqrt{2}k^2}{2}\,\frac{2b^2 + c^2\,s_+}{b^2\,(-b^2 + 4c^2\,s_+)}\,r^{-2}+  O(r^{-4}), \\
v &= -\frac{s_+}{\sqrt{6}} - \frac{\sqrt{6}k^2}{2}\frac{-b^2 + c^2\,s_+}{b^2\,(-b^2 + 4c^2\,s_+)}\,r^{-2} + O(r^{-4}) .
\label{ab2}
\end{align}
\end{lemma}

The proof of this result uses the following lemma:

\begin{lemma}\label{Lemma:MPX}
Let $B_R \subset \RR^n$ with $0<R<\infty$. Assume for some constant $C > 1$ that 
\[
\frac{1}{C} \leq {h}(x) \leq C \text{ in } \RR^n \setminus B_R.
\]
If $u \in C^2(\RR^n \setminus B_R)$ satisfies
\[
-\Delta u + h(x)\,u = O(|x|^{-\alpha})
\]
for some $\alpha > 0$ and if $ u(x) \rightarrow 0$ as $|x| \rightarrow \infty$, then $u = O(|x|^{-\alpha})$, where the big`` $O$'' notation is meant for large $|x|$.
\end{lemma}

\begin{proof}
Let $L = -\Delta + h(x)$. We have 
\[
L(|x|^{-\alpha}) = \alpha(\alpha - n + 2)|x|^{-\alpha-2} + h(x)\,|x|^{-\alpha}.
\]
Hence, by our assumption on $h(x)$, 
\[
\frac{1}{C}|x|^{-\alpha} \leq L(|x|^{-\alpha}) \leq C\,|x|^{-\alpha} \text{ in } \RR^n \setminus B_{2R}.
\]
It thus follows that, there is some large radius $R' > 2R$ and some $C_1 > 0$ such that
\[
L(u - C_1|x|^{-\alpha}) \leq 0 \leq L(u + C_1\,|x|^{-\alpha}) \text{ in } \RR^n \setminus \bar B_{R'}.
\]
Replacing $C_1$ by a larger constant if necessary, we can also assume that
\[
u - C_1|x|^{-\alpha} \leq 0 \leq u + C_1|x|^{-\alpha} \text{ on } \partial B_{R'}.
\]
The assertion follows from the maximum principle.
\end{proof}

\begin{proof}[Proof of Lemma \ref{lemma:asym}.] Let $\hat u  = u - u(\infty)$ and $\hat v = v - v(\infty)$. We have
\begin{align}
\hat u'' + \frac{1}{r}\hat u'
	&= \Big(c_1 + O(|\hat u| + |\hat v| + r^{-2}) \Big)\,\hat u+  \Big(c_2 + O(|\hat u| + |\hat v|)  \Big)\,\hat v  + \frac{k^2\,s_+}{\sqrt{2}r^2},\label{Eq:hw1}\\
\hat v'' + \frac{1}{r}\hat v' 
	&=  \Big(c_2 + O(|\hat u| + |\hat v|)  \Big)\hat u
		+ \Big(c_3 + O(|\hat u| + |\hat v|) \Big) \hat v.\label{Eq:hw0}
\end{align}
where $c_1 = c^2\,s_+^2$, $c_2 = -\frac{\sqrt{3}}{3}(c^2\,s_+^2 - b^2\,s_+)$, and $c_3 = \frac{1}{3}(2b^2\,s_+ + c^2\,s_+^2)$.

Introducing $X = \hat u + \sqrt{3}\hat v $ and $Y =  \sqrt{3}\hat u -\hat v$, we obtain
\begin{align*}
X'' + \frac{1}{r} X' 
	&= \Big(b^2\,s_+ + O(|\hat u| + |\hat v| + r^{-2}) \Big) X
		+ \frac{k^2s_+}{\sqrt{2}r^2} ,\\
Y'' + \frac{1}{r} Y'
	&= \frac{1}{3}\Big(4c^2\,s_+^2 - b^2\,s_+ + O(|\hat u| + |\hat v|+ r^{-2}) \Big)\,Y + \frac{\sqrt{3}k^2\,s_+}{\sqrt{2}r^2}.
\end{align*}
Since both $b^2\,s_+$ and $4c^2\,s_+^2 - b^2\,s_+$ are positive and since $O(|\hat u| + |\hat v|)  = o(1)$ as $r\rightarrow \infty$, Lemma \ref{Lemma:MPX} implies that $|X| \leq C\,r^{-2}$ and $|Y| \leq C\,r^{-2}$. It follows that the above equations of $X$ and $Y$ can be rewritten as 
\begin{align*}
X'' + \frac{1}{r} X' 
	&= b^2\,s_+ X
		+ \frac{k^2s_+}{\sqrt{2}r^2} + O(r^{-4}),\\
Y'' + \frac{1}{r} Y'
	&= \frac{1}{3}(4c^2\,s_+^2 - b^2\,s_+)\,Y + \frac{\sqrt{3}k^2\,s_+}{\sqrt{2}r^2} + O(r^{-4}).
\end{align*}
Thus, the functions $\bar X = X + \frac{k^2}{\sqrt{2}b^2}\,r^{-2}$ and $\bar Y = Y + \frac{3\sqrt{3}k^2}{\sqrt{2}(4c^2\,s_+ - b^2)}r^{-2}$ satisfy
\begin{align*}
\bar X'' + \frac{1}{r} \bar X' 
	&= b^2\,s_+ \bar X
		 + O(r^{-4}),\\
\bar Y'' + \frac{1}{r} \bar Y'
	&= \frac{1}{3}(4c^2\,s_+^2 - b^2\,s_+)\,\bar Y  + O(r^{-4}).
\end{align*}
Again, Lemma \ref{Lemma:MPX} implies that $|\bar X| + |\bar Y| \leq C\,r^{-4}$. Returning to the variables $u$ and $v$, we obtain the desired asymptotic expansion.
\end{proof}

%---------------------------------------------------%
\section{Instability of $k$-radially symmetric solutions}
\label{sec:inst}

In this section we prove the instability of radially $k$-symmetric solutions of \eqref{eq:EL} on the whole $\RR^2$ for $|k|>1$. Note that for any $Q \in H_{loc}^1(\RR^2,\mcS_0)$ satisfying \eqref{BC1}, one has $\mcF(Q) = \infty$. We thus adopt a second variation at $Q$ in a local sense as defined in \eqref{Eq:CLDef}.

%Thus  the second variation  of the functional $\mathcal{\mcF}$ at the point $Y_k$ in the direction $P$ is defined as:
%\bea \label{SV}
%\CL[Y_k](P)&=\frac{1}{2}\frac{d^2}{dt^2}|_{t=0} \int_{\RR^2}\left[\frac{1}{2}|\nabla (Y_k+tP)|^2+f_B(Y_k+tP)-\frac{1}{2}|\nabla Y_k|^2-f_B(Y_k)\right]\,dx
%\non\\
%&=\int_{\RR^2}\frac{1}{2}|\nabla P|^2-\frac{a^2}{2}|P|^2-b^2\tr(P^2Y_k)+\frac{c^2}{2}\left(|Y_k|^2|P|^2+2|\tr(Y_kP)|^2\right)\,dx,
%\eea 
%where $P \in H^1(\RR^2, \mcS_0)$.

%\begin{proposition} \label{prop:instability}\marginnote{why making a new statement? rather than the Theorem?}\marginnote{what is $Y_k$?}
%Let $a^2, b^2, c^2 >0$ be fixed constants and $k \neq 0, \pm 1$. There exists $P \in C_c^\infty(\RR^2, \mcS_0)$ such that the second variation $\CL[Y_k](P) <0$.
%\end{proposition}
\begin{proof}[Proof of Theorem~\ref{thm:main}.]
We follow the ideas from \cite{INSZ_CRAS, INSZ3}. For $|k|>1$, let $Q$ be a $k$-radially symmetric solution of \eqref{eq:EL} on $\RR^2$ subjected to \eqref{BC1}. Then $Q$ has the form \eqref{anY} with $(u,v)$ satisfying \eqref{ODEsystem}-\eqref{bdrycond}.

Let $\eps>0$ be a small parameter. Since $u(\infty) = \frac{s_+}{\sqrt{2}}$, there exists $R > 0$ such that 
\be \label{choic_u}
(1-\eps) u(\infty) \leq u\leq (1+\eps) u(\infty) \quad \textrm{ in } \, \, (R,\infty).
\ee
We take 
$$
P= w(r) h(\varphi) \frac{1}{\sqrt{2}}(n \otimes e_3+e_3\otimes n)
$$ where $n$ is as defined in \eqref{def:n}, $w \in C_c^\infty(R,\infty)$ and
\[
h(\varphi) = \left\{\begin{array}{ll}
	\sin(\frac{\varphi}{2}) & \text{ if $k$ is odd,}\\
	\frac{1}{\sqrt{2}} & \text{ if $k$ is even.}
\end{array}\right.
\]
Then $P \in C_c^\infty(\RR^2,\mcS_0)$.  

%. \marginnote{Arghir: this regularity of $P$ holds only for $k$ odd, but not for $k$ even, because of discontinuity at $\phi=0$. Is there a perturbation that works for both $k$ even and odd?}.

%Note that $P \in H^1(\RR^2, \mcS_0)$  Moreover we have
%In order to systematically study the sign of the second variation we use the basis \eqref{basis}. We define
%$$
%P = \sum_{i=0}^4 w_i \, E_i
%$$
%and plug it into the expression for the second variation to obtain
%\begin{align} 
%2 {\CL}_R[Y_k](P)&=\int_{B_R} \sum_{i=0}^4 |\partial_r w_i |^2  \non \\ 
%&+\frac{1}{r^2} \left( |\partial_\phi w_0|^2+ |\partial_\phi w_1 - k w_2|^2 + 
%|\partial_\phi w_2 + k w_1|^2 + |\partial_\phi w_3 - \frac{k}{2} w_4|^2 + |\partial_\phi w_4 + \frac{k}{2} w_3|^2\right) \non\\
%&+ \left(-{a^2} + c^2 (u^2+v^2) \right)\sum_{i=0}^4 |w_i |^2+2 {c^2}\left( v w_0 + u w_1 \right)^2 \\
% &-b^2 \left( \frac{\sqrt{2}}{2} u (w_3^2 -w_4^2) + \frac{1}{ \sqrt{6}} v ( w_3^2 + w_4^2 +2 w_0^2 - 2w_1^2 -2 w_2^2) -\frac{4}{\sqrt{6}} u w_0 w_1 \right),
%\end{align}
%Now we assume $w_0=w_1=w_2=w_4=0$ and $w_3 = w(r)$ to obtain
We have  

$$\frac{\partial P}{\partial r}=w'h\frac{1}{\sqrt{2}}\left(n\otimes e_3+e_3\otimes n\right),\quad \frac{\partial P}{\partial \varphi}=wh'\frac{1}{\sqrt{2}}\left(n\otimes e_3+e_3\otimes n\right)+wh\frac{k}{2\sqrt{2}}\left(m\otimes e_3+e_3\otimes m\right)$$ so 

$$\frac{|\nabla P|^2}{2}=\frac{|\partial_r P|^2}{2}+\frac{|\partial_\varphi P|^2}{2r^2}=\frac{(w'h)^2}{2}++\frac{w^2 (h')^2}{2r^2}+\frac{w^2k^2h^2}{8r^2}.$$

We also have:

$$\int_0^{2\pi} h^2(\varphi)\,d\varphi=\pi, \int_0^{2\pi} (h'(\varphi))^2\,d\varphi=c_k=\frac{1}{2}(1+(-1)^{k+1}),$$

$$|P^2|=w^2h^2, |Q|^2=u^2+v^2, \tr(PQ)=0$$

$$P^2=\frac{w^2h^2}{2}(n\otimes n+e_3\otimes e_3), \tr(P^2Q)=\frac{w^2h^2}{2}(\frac{u}{\sqrt{2}}+\frac{v}{\sqrt{6}}),$$

hence \eqref{Eq:CLDef} becomes:
\begin{align*}
\frac{2}{\pi} {\CL} [Q](P)
	& =  \int_0^\infty \Big\{ |w'|^2 +  \frac{k^2+ c_k }{4 r^2} w^2 + \left(-a^2 +\frac{2}{\sqrt{6}} b^2 v + c^2 (u^2 +v^2) \right) w^2\\
	&\qquad\qquad -   \frac{b^2}{\sqrt{2}} \left( u+ \sqrt{3} v \right) w^2  \Big\} r\, dr,
\end{align*}
where 
\[
c_k = \frac{1}{2}(1 + (-1)^{k+1}).
\]
We now use the Hardy decomposition trick as in \cite{INSZ_CRAS, INSZ3} by setting $w=u \xi$ with $\xi \in C_c^\infty(R,\infty)$. Then:

\begin{align*}
\frac{2}{\pi} {\CL} [Q](P)& =\int_0^\infty \Big\{ |u'\xi|^2 +|u\xi'|^2+2uu'\xi\xi'+\frac{k^2+ c_k }{4 r^2} u^2\xi^2 + (u''+\frac{u'}{r}-\frac{4k^2u^2}{4r^2})u\xi^2\\
	&\qquad\qquad -   \frac{b^2}{\sqrt{2}} \left( u+ \sqrt{3} v \right) u^2\xi^2  \Big\} r\, dr,\\
 &=\int_0^\infty \left\{  |\xi'|^2  - \frac{3 k^2- c_k}{4r^2} \xi^2 -   \frac{b^2}{\sqrt{2}} \left( u+ \sqrt{3} v \right) \xi^2 \right\} u^2 r\, dr .
\end{align*} where for the first equality we used the equation \eqref{ODEsystem} for $u$ and for the second equality we integrated by parts the term $\int_0^\infty u'' u\xi^2 r\,dr$.
From Lemma~\ref{lemma:asym} we know
$$
u+ \sqrt{3} v = -\frac{k^2}{\sqrt{2} b^2 r^2}  + O(r^{-4}) \quad \hbox{ as }r \to \infty.
$$  
Therefore, by replacing $R$ by a larger constant if necessary, we can assume that
$$ \frac{b^2}{\sqrt{2}} \left( u+ \sqrt{3} v \right)\geq \frac{-11k^2+c_k}{20r^2}  \quad \textrm{ in } \, \, (R,\infty).
$$
Hence, for any $\xi \in  C_c^\infty (R,\infty)$, we deduce by \eqref{choic_u} and $k^2-c_k\geq 3$:
$$
\frac{2}{\pi} {\CL} [Y_k](P) \leq \int_R^\infty \left\{  |\xi'|^2  - \frac{ k^2- c_k }{5r^2} \xi^2 \right\}  u^2 r\, dr \leq {(1+\eps)^2 u(\infty)^2} \int_R^\infty \left\{  |\xi'|^2  - \frac{ 1}{2r^2} \xi^2 \right\}  r\, dr .
$$
It is not difficult to find a test function $\xi_0 \in C_c^\infty (R,\infty)$ such that \footnote{For example, take $\xi_0$ to be a smoothing of $\sin(\frac{\ln r}{2})\,\mathbf{1}_{(\exp(2n\pi),\exp(2(n+1)\pi))}$ for some $n$ sufficiently large.}
$$
\int_R^\infty \left\{  |\xi_0'|^2  - \frac{ 1}{2r^2} \xi_0^2 \right\}  r\, dr <0.
$$
The result follows immediately.
\end{proof}

\begin{remark} Theorem~\ref{thm:main} and its proof provide an insight into the stability of the $k$-radially solution on finite domains $B_R(0)$ for $R$ small, respectively $R$ large.

\begin{itemize}
\item {\it Case 1: if $R$ is small (ie, $R\leq R_0(a^2, b^2, c^2)$)} then one can use the Poincar\'e inequality to show that  the solution of the PDE system \eqref{eq:EL} with boundary conditions \eqref{BC1f} is unique  (see  for instance in the related  Ginzburg-Landau framework Thm. $VIII.7$, p. $98$,\cite{vortices}). This unique solution must necessarily be the global minimizer of $\mcF$ on $\Omega=B_R$ subject to \eqref{BC1f}.
 Therefore, it coincides with the $k$-radially solution of Proposition 2.3 (write \ref{pro:sol_rad_sym}) and it is stable as global minimizer of $\mcF$. 

\item {\it Case 2: if $R$ is large (ie, $R\geq R_1(a^2, b^2, c^2)$)} then for $|k|>1$ the  $k$-radially solutions obtained in Proposition~\ref{prop:positivlocalmin} are expected to be unstable, since the  solutions $(u_R,v_R)$ of the ODE system   on the finite domain $[0,R]$ will suitably approximate the solution in the whole space obtained in Proposition~\ref{prop:inf}.
\end{itemize}
\end{remark}

%We are ready to conclude with the proof of Theorem~\ref{thm:main}.
%\begin{proof}
%Let $a^2, b^2, c^2 >0$, $k \in \ZZ \setminus \{ 0, \pm 1\}$ and $Q$ be $k$-radially symmetric solution of the Euler-Lagrange equation \eqref{eq:EL}. Combining Propositions~\ref{2ode}, \ref{prop:inf} and \ref{prop:instability} we obtain the result.
%\end{proof}
%

\appendix
%---------------------------------------------------%

\section{Remarks on the case $b^2 = 0$}
\label{app:b0}

In this appendix, we collect some known results regarding the case $b^2 = 0$ and $k \in \ZZ \setminus \{0\}$.
\smallskip
 In \cite{DRSZ}, it was shown that on a finite disk the system \eqref{ODEsystem} and \eqref{BCodeFinite} has a unique solution $(u,v)$ with the sign invariance $u > 0$ and $v < 0$. Furthermore, $Q = u\,E_1 + v\,E_0$ is the unique global minimizer of the full Landau-de Gennes energy $\mcF$ subjected to the boundary condition \eqref{BC1f}.

\smallskip
 For infinite domain, the situation is different. We have:

\begin{theorem}
Assume that $b^2 = 0$, $a^2 > 0$, $c^2 > 0$ and $k \neq 0$. There is no solution of the boundary value problem \eqref{ODEsystem}-\eqref{bdrycond} which satisfies $u > 0$ and $v < 0$ in $(0,\infty)$. 
\end{theorem}

\begin{proof}
Indeed, assume by contradiction that there exists a solution $(u,v)$ of \eqref{ODEsystem} on $(0,\infty)$ subject to \eqref{bdrycond} with $v < 0$ in $(0,\infty)$. By \eqref{bds:uvmax} and \eqref{def:s_+}, we have that $-a^2+c^2(u^2+v^2)\le 0$ in $(0, \infty)$. Hence, the equation \eqref{ODEsystem} for $v\leq 0$ implies $(rv')'\ge 0$ for every $r>0$. Since $v'(0)=0$, we deduce $r\mapsto rv'(r)$ is a nonnegative and nondecreasing function. It follows that, for any $s > r > 0$, $v(s) \geq v(r) + rv'(r)\,\ln \frac{s}{r}$. Fixing $r$ and taking a limit  $s \to \infty$ it is clear that starting from some point $s_0$ function $v(s)$ becomes positive. Since $v$ is negative in $(0,\infty)$, this implies that $v' \equiv 0$ in $(0,\infty)$. By \eqref{bdrycond}, we thus have $v\equiv-\frac{s_+}{\sqrt{6}}$ in $(0, \infty)$. Using the second equation in \eqref{ODEsystem}, we obtain $-a^2+c^2(u^2+v^2)=0$, and so, by \eqref{bdrycond}, $u \equiv  \frac{s_+}{\sqrt{2}}$ is constant. This contradicts the first equation in \eqref{ODEsystem}. 
\end{proof}

\section*{Acknowledgment.} The authors gratefully acknowledge the hospitality and partial support of the Mathematisches Forschungsinstitut Oberwolfach, Centre International de Rencontres Math\'{e}matiques, and Institut Henri Poincar\'{e}, where parts of this work were carried out. V.S. and A.Z. thank Dr. Jonathan Robbins for insightful discussions. R.I. acknowledges partial support by the ANR project ANR-14-CE25-0009-01, V.S. acknowledges partial support by EPSRC grant EP/K02390X/1.

\bibliographystyle{acm}
\bibliography{paris,LiquidCrystals}

\end{document}